# On mean central limit theorems for stationary sequences

Jérôme Dedecker[a] and Emmanuel Rio[b]

[a] *Laboratoire de Statistique Théorique et Appliquée, Université Paris 6, 175 rue du Chevaleret, 75013 Paris, France. E-mail: dedecker@ccr.jussieu.fr*
[b] *Laboratoire de mathématiques, UMR 8100 CNRS, Bâtiment Fermat, 45 Avenue des Etats-Unis, 78035 VERSAILLES Cedex, France. E-mail: rio@math.uvsq.fr*

**Abstract.** In this paper, we give estimates of the minimal $\mathbb{L}^1$ distance between the distribution of the normalized partial sum and the limiting Gaussian distribution for stationary sequences satisfying projective criteria in the style of Gordin or weak dependence conditions.

**Résumé.** Dans cet article, nous donnons des majorations de la distance minimale $\mathbb{L}^1$ entre la loi de la somme normalisée et sa loi limite gaussienne pour des suites stationnaires satisfaisant des critères projectifs à la Gordin ou des conditions de dépendance faible.

*MSC:* 60F05

*Keywords:* Mean central limit theorem; Wasserstein distance; Minimal distance; Martingale difference sequences; Strong mixing; Stationary sequences; Weak dependence; Rates of convergence; Projective criteria

## 1. Introduction

Let $X_1, X_2, \dots,$ be a sequence of real-valued random variables (r.v.) with mean zero and finite variance. Let $S_n = X_1 + X_2 + \cdots + X_n$. By $F_n$ we denote the distribution function (d.f.) of $n^{-1/2}S_n$. Let $\Phi_\sigma$ be the d.f. of the $\mathcal{N}(0,\sigma^2)$-distribution. For independent and identically distributed (i.i.d.) r.v.'s, according to the central limit theorem (CLT), $F_n(x)$ converges to $\Phi_\sigma(x)$ uniformly for $x$ in $\mathbb{R}$, where $\sigma$ is the standard deviation of $X_1$. Agnew [1] proved that the convergence also holds in $L^r(\mathbb{R})$ for $r > 1/2$. Agnew's result is called mean CLT in the case $r = 1$. Let then $\rho_n^{(r)} = \|F_n - \Phi_\sigma\|_r$. For $r = 1$ and $r = 2$ and i.i.d. random variables with finite absolute third moment, Esseen [11] proved that $n^{1/2}\rho_n^{(r)}$ converges to some explicit constant $A_r(F)$ depending only on the distribution function $F$ of $X_1$ (Theorems 3.2 and 4.2 in [11]). In particular, Esseen's results imply that

$$\rho_n^{(r)} = \mathrm{O}(n^{-1/2}) \quad \text{as } n \to \infty. \tag{1.1}$$

Next Zolotarev [29] obtained the upper bound $A_1(F) \le \mathbb{E}(|X_1|^3)/(2\sigma^2)$. The proofs of these results are based on the method of characteristic functions (cf. [18] for more details).

---

The case $r = 1$ is of special interest, since $\rho_n^{(1)}$ is exactly the minimal distance between $n^{-1/2}S_n$ and a r.v. with distribution $\mathcal{N}(0, \sigma^2)$ in $\mathbb{L}^1$ (cf. [10], Section 11.8, Problems 1 and 2). Now let

$$d_1(X, Y) = \sup_{f \in \Lambda_1(\mathbb{R})} \mathbb{E}(f(X) - f(Y)), \tag{1.2}$$

where $\Lambda_1(\mathbb{R})$ is the set of 1-Lipschitzian functions from $\mathbb{R}$ to $\mathbb{R}$. Applying the Kantorovich–Rubinstein theorem we also have that $\rho_n^{(1)} = d_1(n^{-1/2}S_n, \sigma Y)$ if $Y$ is a $\mathcal{N}(0, 1)$-distributed random variable.

In this paper we are interested in extensions of (1.1) for $r = 1$ to sequences of dependent random variables. This subject was studied by Sunklodas [28] in the case of uniformly mixing (in the sense of Ibragimov) stationary sequences of real-valued random variables. Using the Stein method, he reached the rate of convergence $\mathrm{O}(n^{-1/2}(\log n)^2)$ in (1.1) for geometrically mixing sequences of random variables with finite eight moments. A different approach to get rates of convergence in the CLT is Bergström's [2] inductive proof of the Berry–Esseen theorem, based on the Lindeberg method. Starting from Bergström's recursion argument, Bolthausen [4] obtained exact rates of uniform convergence for martingale difference arrays. Rio [25] adapted Bergström's method to weakly dependent sequences and obtained the Berry–Esseen theorem for stationary and uniformly bounded sequences of real-valued r.v.'s satisfying the condition $\sum_k k\varphi(k) < \infty$, where $(\varphi(k))_k$ denotes the sequence of uniform mixing coefficients of the sequence $(X_i)_{i \in \mathbb{N}}$, in the sense of Ibragimov (confer [18] for an exact definition of these coefficients). This result was extended to the multivariate case by Jan [19], Theorem 9. Jan also weakened the notion of weak dependence involved in Rio's paper (cf. Theorem 1 in [20] for more details). However the dependence coefficients in [19] are too restrictive for the applications to some dynamical systems, such as Sinai's billiard. Pène [22] noticed that the inductive proof of Jan [19] can be adapted to get the rate of convergence $\mathrm{O}(n^{-1/2})$ for the minimal $\mathbb{L}^1$-distance in the multivariate CLT for stationary sequences satisfying some dependence conditions. In particular her result applies to sums of bounded r.v.'s defined from dynamical systems (such as Sinai's billiard) or strongly mixing sequences in the sense of Rosenblatt. For example, Pène's result yields (1.1) (with $r = 1$) for stationary sequences of bounded random variables $(X_i)_{i \in \mathbb{N}}$ satisfying the condition $\sum_k k\alpha(k) < \infty$, where $(\alpha(k))_k$ denotes the sequence of strong mixing coefficients of $(X_i)_{i \in \mathbb{Z}}$ in the sense of Rosenblatt (confer [18] for a definition of these coefficients).

We now describe the contents of our paper. Our aim is to provide rates of convergence in the mean CLT for stationary sequences of real-valued r.v.'s satisfying either projective criteria in the style of Gordin [14] or weak dependence conditions.

In Section 2, we give bounds in the stationary case involving $\mathbb{L}^p$-norms of conditional expectations. Let $(X_i)_{i \in \mathbb{Z}}$ be a stationary sequence of real-valued random variables, $\mathcal{M}_k = \sigma(X_i\colon\ i \le k)$ and $E_k$ denote the conditional expectation with respect to $\mathcal{M}_k$. In Section 2.1, we obtain in Theorem 2.1 the rate of convergence $\mathrm{O}(n^{-1/2} \log n)$ in the mean CLT for stationary and ergodic martingale differences sequences $(X_i)_{i \in \mathbb{Z}}$ with finite absolute third moments satisfying the projective conditions

$$\sup_{m>0} \left\| \sum_{k=1}^{m} E_0(X_k^2 - \sigma^2) \right\|_1 < \infty \quad \text{and} \quad \sup_{m>0} \left\| \sum_{k=1}^{m} X_0 E_0(X_k^2 - \sigma^2) \right\|_1 < \infty, \tag{1.3}$$

where $\sigma^2 = \operatorname{Var} X_0$. In Section 2.2, we generalize Theorem 2.1 to ergodic stationary sequences satisfying projective criteria. In Section 2.3 we give some applications to bounded sequences. For example, assuming that the series $\sum_{k>0} E_0(X_k)$ converges in $\mathbb{L}^1$, Theorem 2.3 provides rates of convergence in the mean CLT as soon as $E_0(S_m^2/m)$ converges to $\sigma^2$ in $\mathbb{L}^1$. This condition appears in the conditional CLT of Dedecker and Merlevède [6] and is rather mild. For example the rate of convergence $\mathrm{O}(n^{-1/2} \log n)$ is obtained under the projective conditions

$$\sum_{m>0} \left\| \sum_{k \ge m} E_0(X_k) \right\|_1 < \infty \quad \text{and} \quad \sup_{m>0} \|E_0(S_m^2 - m\sigma^2)\|_1 < \infty. \tag{1.4}$$

Again the proofs are based on the Lindeberg method at order three.

In Section 3, we give projective conditions or weak dependence conditions implying (1.1) for $r = 1$. Conditions (1.3) and (1.4) involve conditional second moments. It seems difficult to get the optimal rate of convergence $\mathrm{O}(n^{-1/2})$ under second-order conditions (at least for the Berry–Esseen theorem: cf. [25] and [4], Theorem 4). Therefore our results hold under projective conditions on the monoms of degree three. For example, (1.1) holds for stationary bounded martingale difference sequences under the projective conditions

$$\sum_{k>0} \|E_0(X_k^2) - \sigma^2\|_1 < \infty \quad \text{and} \quad \sum_{k>0} \sup_{j\geq k} \|E_0(X_k X_j^2) - \mathbb{E}(X_k X_j^2)\|_1 < \infty. \tag{1.5}$$

For stationary sequences, one needs to strengthen (1.5): we obtain (1.1) for stationary sequences of bounded r.v.'s under the projective conditions

$$\sum_{k>0} k \sup_{i\geq j\geq k} \|E_0(X_k X_j^\alpha X_i^\beta) - \mathbb{E}(X_k X_j^\alpha X_i^\beta)\|_1 < \infty, \quad \text{with } (\alpha,\beta) \in \{0,1\}^2, \tag{1.6}$$

which can also be deduced from Theorem 1.1 in [22]. It is worth noticing that the Berry–Esseen type Theorem 9 in [19] requires $\mathbb{L}^\infty$-norms instead of $\mathbb{L}^1$-norms in (1.6). The proofs of these results are based on the Lindeberg method at order four. Therefore, in the unbounded case, the results hold for sequences of random variables with finite fourth moments (cf. Theorems 3.1(a) and 3.2(a) for detailed conditions). For example, Theorem 3.1(a) applied to strongly mixing and stationary sequences yields the rate of convergence $\mathrm{O}(n^{-1/2})$ in the mean CLT if there exists some $p > 1$ such that

$$\sum_{k>0} k^{(p+1)/(p-1)} \alpha(k) < \infty \quad \text{and} \quad \mathbb{E}(|X_0|^{ap}) < \infty, \tag{1.7}$$

where $a = 4$. By contrast, the Berry–Esseen type theorem for functionals of stationary discrete Markov chains due to Bolthausen [3] holds under condition (1.7) with $a = 3$. In order to improve Theorems 3.1(a) and 3.2(a) in the case of strongly mixing sequences we adapt the truncation method in [24] to our context. We then get the rate $\mathrm{O}(n^{-1/2})$ in the mean CLT under the strong mixing condition

$$\sum_{k>0} k^b \int_0^{\alpha(k)} Q_{|X_0|}^3(u)\,\mathrm{d}u < \infty, \tag{1.8}$$

where $Q_{|X_0|}$ denotes the quantile function of $|X_0|$ and $b = 1$. This condition is implied by (1.7) with $a = 3$, so that our result holds under Bolthausen's [3] condition. Moreover, for stationary strongly mixing martingale difference sequences, we prove that (1.1) holds for $p = 1$ under condition (1.8) with $b = 0$. In Section 5 we give two classical examples of non irreducible Markov chains to which our results apply.

## 2. Projective criteria for stationary sequences

Throughout the paper, $Y$ is a $\mathcal{N}(0,1)$-distributed random variable.

We shall use the following notations. Let $(\Omega, \mathcal{A}, \mathbb{P})$ be a probability space, and $T : \Omega \mapsto \Omega$ be a bijective bimeasurable transformation preserving the probability $\mathbb{P}$. An element $A$ is said to be invariant if $T(A) = A$. We denote by $\mathcal{I}$ the $\sigma$-algebra of all invariant sets. Let $\mathcal{M}_0$ be a sub-$\sigma$-algebra of $\mathcal{A}$ satisfying $\mathcal{M}_0 \subseteq T^{-1}(\mathcal{M}_0)$ and define the nondecreasing filtration $(\mathcal{M}_i)_{i\in\mathbb{Z}}$ by $\mathcal{M}_i = T^{-i}(\mathcal{M}_0)$. Let $\mathcal{M}_\infty = \bigvee_{i\in\mathbb{Z}} \mathcal{M}_i$. Denote by $E_i$ the conditional expectation with respect to $\mathcal{M}_i$.

Let $X_0$ be a $\mathcal{M}_0$-measurable and centered random variable. Throughout the sequel, the sequence $\mathbf{X} = (X_i)_{i\in\mathbb{Z}}$ is defined by $X_i = X_0 \circ T^i$. From the definition the sequence $(X_i)_{i\in\mathbb{Z}}$ is adapted to the filtration $(\mathcal{M}_i)_{i\in\mathbb{Z}}$.

### 2.1. Martingale difference sequences

In this section we obtain rates of convergence of the order of $n^{-1/2}\log n$ in the mean CLT for stationary martingale difference sequences. In order to obtain these rates of convergence, we will just need a projective condition on the variables $X_l^2$, as in [19]. We first recall Jan's results concerning the rates of convergence for the uniform distance between the distribution functions.

Assume that $(X_i)_{i\in\mathbb{Z}}$ is a stationary martingale difference sequence in $\mathbb{L}^3$ such that $\mathbb{E}(X_0^2)=\sigma^2$ and

$$\sum_{l>0}\|E_0(X_l^2-\sigma^2)\|_{3/2}<\infty. \tag{2.1}$$

Then, by Theorem 6 in [19], if $Y$ is $\mathcal{N}(0,1)$-distributed,

$$\sup_{t\in\mathbb{R}}|\mathbb{P}(n^{-1/2}S_n\le t)-\mathbb{P}(\sigma Y\le t)|=\mathrm{O}(n^{-1/4}). \tag{2.2}$$

Under projective conditions related to (2.1), the rate of convergence in the mean central limit theorem is at least $\mathrm{O}(n^{-1/2}\log n)$ as shown in Theorem 2.1 below.

**Theorem 2.1.** *Let $(X_i)_{i\in\mathbb{Z}}$ be a stationary martingale difference sequence in $\mathbb{L}^3$, such that $\mathbb{E}(X_0^2|\mathcal{I})=\mathbb{E}(X_0^2)=\sigma^2$ almost surely. Let $\Lambda=\sigma^{-2}\mathbb{E}|X_0|^3$ and $U_m=E_0(X_1^2+\cdots+X_m^2)-m\sigma^2$. Then*

(a) $d_1(S_n,\sigma\sqrt{n}Y)\le\frac{13\sigma}{6}+\frac{\Lambda}{6}\log(1+2n)+\sum_{m=1}^{[\sqrt{2n}]}\frac{\|X_0U_m\|_1+2\sigma\|U_m\|_1}{m\sigma^2}$.
(b) *If* $\sup_{m>0}(\|X_0U_m\|_1+\|U_m\|_1)<\infty$, *then* $d_1(S_n,\sigma Y\sqrt{n})=\mathrm{O}(\log n)$.

***Remark 2.1.*** *From the ergodic theorem, $(U_m/m)$ converges a.s. and in $\mathbb{L}^1$ to 0 as $m$ tends to $\infty$. Since $X_0\in\mathbb{L}^3$, it follows that the sequence $(X_0U_m/m)_m$ is uniformly integrable. Hence, under the assumptions of Theorem 2.1,*

$$\lim_{m\to\infty}m^{-1}(\|X_0U_m\|_1+\|U_m\|_1)=0.$$

*Therefore Theorem 2.1*(a) *provides a rate of convergence in the mean CLT. For example, if $\|X_0E_0(X_l^2-\sigma^2)\|_1=\mathrm{O}(l^{-\delta})$ and $\|E_0(X_l^2-\sigma^2)\|_1=\mathrm{O}(l^{-\delta})$ for some $\delta$ in $]0,1[$, then the rate of convergence in the mean CLT is of the order of $n^{-\delta/2}$. If Jan's condition (2.1) holds, then* (b) *yields the rate of convergence $\mathrm{O}(n^{-1/2}\log n)$ in the mean CLT. For bounded random variables* (b) *holds as soon as the series $\sum_{l>0}E_0(X_l^2-\sigma^2)$ converges in $\mathbb{L}^1$.*

**Proof of Theorem 2.1.** We prove Theorem 2.1 in the case $\sigma=1$. The general case follows by dividing the r.v.'s by $\sigma$.

Let $(Y_i)_{i\in\mathbb{N}}$ be a sequence of independent random variables with normal distribution $\mathcal{N}(0,1)$. Suppose furthermore that the sequence $(Y_i)_{i\in\mathbb{N}}$ is independent of $(X_i)_{i\in\mathbb{N}}$. Let $Y$ be a $\mathcal{N}(0,1)$-distributed random variable, independent of the above defined sequences. Let $T_n=Y_1+Y_2+\cdots+Y_n$. For any 1-Lipschitzian function $f$, let $\Delta(f)=\mathbb{E}(f(S_n)-f(T_n))$. From (1.2), we have to bound $\Delta(f)$. Clearly

$$\Delta(f)=\mathbb{E}(f(S_n)-f(T_n))\le\mathbb{E}(f(S_n+Y)-f(T_n+Y))+2\mathbb{E}|Y|. \tag{2.3}$$

In order to bound up the term on right-hand side, we apply the Lindeberg method.

***Notation 2.1.*** *Set $f_k(x)=\mathbb{E}(f(x+Y+T_n-T_k))$. Let $S_0=0$, and, for $k>0$, let $\Delta_k=f_k(S_{k-1}+X_k)-f_k(S_{k-1}+Y_k)$.*

Since the sequence $(Y_i)_{i\in\mathbb{N}}$ is independent of the sequence $(X_i)_{i\in\mathbb{N}}$,

$$\mathbb{E}(f(S_n+Y)-f(T_n+Y))=\sum_{k=1}^{n}\mathbb{E}(\Delta_k). \tag{2.4}$$

Next the functions $f_k$ are $C^\infty$. Consequently, from the Taylor integral formula at orders three and four,

$$\Delta_k = f_k'(S_{k-1})(X_k - Y_k) + \frac{1}{2} f_k''(S_{k-1})(X_k^2 - Y_k^2) - \frac{1}{6} f_k^{(3)}(S_{k-1}) Y_k^3 + R_k,$$

with

$$R_k \le \frac{1}{6} \|f_k^{(3)}\|_\infty |X_k|^3 + \frac{1}{24} \|f_k^{(4)}\|_\infty Y_k^4. \tag{2.5}$$

Consequently, for any 1-Lipschitzian function $f$,

$$\Delta(f) \le 2\mathbb{E}|Y| + \sum_{k=1}^n \mathbb{E}(R_k) + \sum_{k=1}^n \mathbb{E}\biggl( f_k'(S_{k-1}) X_k + \frac{1}{2} f_k''(S_{k-1})(X_k^2 - 1)\biggr). \tag{2.6}$$

The terms $\mathbb{E}(f_k'(S_{k-1})X_k)$ vanish under the martingale assumption. To bound up the other terms appearing in (2.6), we need to bound up the derivatives of $f_k$. This will be done via the lemma below.

**Lemma 2.1.** *Let $f$ be a 1-Lipschitzian function, $Y$ be a standard normal and $B$ be a real-valued random variable, independent of $Y$. Then*

$$\left| \frac{\mathrm{d}^i}{\mathrm{d}x^i} \mathbb{E} f(x + tY + B) \right| \le t^{1-i} \|\phi^{(i-1)}\|_1 \quad \textit{for any } t > 0$$

*and any positive integer $i$, where $\phi$ denotes the density of $Y$.*

**Proof of Lemma 2.1.** Let $\phi_t$ be the density of $tY$. Then

$$\mathbb{E} f(x + tY + B) = \mathbb{E}(f * \phi_t(x + B)).$$

Since $f$ is 1-Lipschitzian, the Stieltjes measure $\mathrm{d}f$ of $f$ is absolutely continuous with respect to the Lebesgue measure $\lambda$ and $f' = \mathrm{d}f/\mathrm{d}\lambda$ belongs to $[-1,1]$. Next $(f * \phi_t)^{(i)} = f' * \phi_t^{(i-1)}$, and consequently

$$\left| \frac{\mathrm{d}^i}{\mathrm{d}x^i} \mathbb{E} f * \phi_t(x + B) \right| \le \|f'\|_\infty \|\phi_t^{(i-1)}\|_1.$$

Since $\phi_t^{(i-1)}(x) = t^{-i} \phi^{(i-1)}(x/t)$, it implies Lemma 2.1. □

Noting that

$$\|\phi'\|_1 = \sqrt{\frac{2}{\pi}} \le \frac{4}{5}, \qquad \|\phi''\|_1 = \sqrt{\frac{8}{\pi \mathrm{e}}} \le 1, \qquad \|\phi^{(3)}\|_1 = \sqrt{\frac{2}{\pi}} + \sqrt{\frac{32}{\pi e^3}} \le \frac{8}{5}, \tag{2.7}$$

and applying Lemma 2.1 with $t = \sqrt{n-k+1}$, we infer from (2.5) that

$$\mathbb{E}(R_k) \le \left(\frac{\Lambda}{6}\right)(n-k+1)^{-1} + \left(\frac{1}{5}\right)(n-k+1)^{-3/2}. \tag{2.8}$$

Summing on $k$, we infer from (2.8) that

$$\sum_{k=1}^n \mathbb{E}(R_k) + 2\mathbb{E}|Y| \le \rho(n) \quad \text{with } \rho(n) = \frac{13}{6} + \frac{\Lambda}{6} \log(1 + 2n). \tag{2.9}$$

The control of the main term in (2.6) is derived from the lemma below.

**Lemma 2.2.** *Let $Z_0$ be an integrable random variable with zero mean. Set $Z_k = Z_0 \circ T^k$ and let $W_m = \sum_{l=1}^m E_0(Z_l)$. Then, for $s=2$ or $s=3$,*

$$\sum_{k=1}^n \mathbb{E}(f_k^{(s)}(S_{k-1})Z_k) \le \sum_{m=1}^{[\sqrt{2n}]} 2m^{1-s}(\|X_0 W_m\|_1 + 2\|W_m\|_1).$$

**Proof of Lemma 2.2.** We first divide $[1,n]$ into blocks of nonincreasing length.

***Notation 2.2.*** *Define the decreasing sequence of integers $(n_i)_{i\ge 0}$ by $n_0 = n$ and $n_i = \max(0, n_{i-1} - i)$ for $i > 0$. Let $p$ be the first integer such that $n_p = 0$. Set $m_i = i$ for $i < p$ and $m_p = n_{p-1}$.*

Next fix $i$ in $[1,p]$. Let then $k$ be any integer in $]n_i, n_{i-1}]$. Writing

$$f_k^{(s)}(S_{k-1}) = f_{n_i+1}^{(s)}(S_{n_i}) + \sum_{j=n_i+1}^{k-1} (f_{j+1}^{(s)}(S_j) - f_j^{(s)}(S_{j-1})), \tag{2.10}$$

we get that

$$\sum_{k=n_i+1}^{n_{i-1}} \mathbb{E}(f_k^{(s)}(S_{k-1})Z_k) = D_i + \sum_{j=n_i+1}^{n_{i-1}-1} D_{i,j}, \tag{2.11}$$

where

$$D_i = \mathbb{E}\left( f_{n_i+1}^{(s)}(S_{n_i}) \sum_{k=n_i+1}^{n_{i-1}} E_{n_i}(Z_k) \right)$$

$$D_{i,j} = \mathbb{E}\left( (f_{j+1}^{(s)}(S_j) - f_j^{(s)}(S_{j-1})) \sum_{k=j+1}^{n_{i-1}} E_j(Z_k) \right).$$

By definition of the sequence $(Z_k)_k$, for any integer $j$ and any positive $m$,

$$E_j(Z_{j+1} + Z_{j+2} + \cdots + Z_{j+m}) = W_m \circ T^j.$$

Hence, from Lemma 2.1 applied with $t = (n-n_i)^{1/2}$ and $B = 0$, for any $i < p$,

$$D_i \le (n - n_i)^{(1-s)/2} \|W_i\|_1. \tag{2.12}$$

Moreover $D_p = 0$ from the centering assumption on the random variables $Z_k$. Now, by definition, $n - n_i = i(i+1)/2$ for $i < p$. Hence, from (2.12),

$$\sum_{i=1}^{p-1} D_i \le 2 \sum_{i=1}^{p-1} i^{1-s} \|W_i\|_1, \quad \text{where } -1 + \sqrt{2n} < p < 1 + \sqrt{2n}. \tag{2.13}$$

Next we bound up $D_{i,j}$. From the elementary equality

$$f_{j+1}^{(s)}(S_j) - f_j^{(s)}(S_{j-1}) = E_j(f_{j+1}^{(s)}(S_{j-1} + X_j) - f_{j+1}^{(s)}(S_{j-1} + Y_j))$$

we get that

$$|f_{j+1}^{(s)}(S_j) - f_j^{(s)}(S_{j-1})| \le \|f_{j+1}^{(s+1)}\|_\infty E_j|X_j - Y_j| \le (n-j)^{-s/2}(|X_j| + 1),$$

whence

$$D_{i,j} \le (n-j)^{-s/2}\mathbb{E}((|X_0|+1)|W_{n_{i-1}-j}|). \tag{2.14}$$

Fix $n_{i-1}-j=m$. Then $m_i > m > 0$ and $2(n-j) = i(i-1)+2m \ge (i-1/2)^2$. Hence, from the above inequality (recall that $m_i = i$ for $i < p$)

$$\sum_{i=1}^{p}\sum_{j=n_i+1}^{n_{i-1}-1} D_{i,j} \le \sum_{m=1}^{p-1}(\|X_0W_m\|_1 + \|W_m\|_1)\sum_{i=m+1}^{p} 2^{s/2}\left(i-\frac{1}{2}\right)^{-s}.$$

Now, from the convexity of $x^{-s}$ on $]0,+\infty[$,

$$\sum_{i=m+1}^{p}\left(i-\frac{1}{2}\right)^{-s} \le \int_m^p x^{-s}\,\mathrm{d}s \le \frac{1}{s-1}m^{1-s}$$

whence

$$\sum_{i=1}^{p}\sum_{j=n_i+1}^{n_{i-1}-1} D_{i,j} \le \sum_{m=1}^{p-1} 2m^{1-s}(\|X_0W_m\|_1 + \|W_m\|_1). \tag{2.15}$$

From (2.11), (2.13) and (2.15), we get Lemma 2.2. □

Theorem 2.1(a) follows from both (2.6), (2.9) and Lemma 2.2 applied to $Z_0 = X_0^2 - 1$. Theorem 2.1(b) is a consequence of (a). □

### 2.2. *Projective criteria*

In this section we give estimates of the rates of convergence in the mean CLT for stationary sequences satisfying projective $\mathbb{L}^1$-criteria in the style of Gordin [14]. Our main result is Theorem 2.2 below.

**Theorem 2.2.** *Let $(X_i)_{i\in\mathbb{Z}}$ be a stationary sequence of centered random variables in $\mathbb{L}^3$ such that $\mathbb{E}(X_0X_k|\mathcal{I}) = \mathbb{E}(X_0X_k)$ a.s. for any integer $k$. Suppose furthermore that the sequence $X_0E_0(S_n)$ converges in $\mathbb{L}^1$. Then the series $\mathbb{E}(X_0^2) + 2\sum_{k=1}^{\infty}\mathbb{E}(X_0X_k)$ is convergent to some nonnegative real $\sigma^2$. Let*

$$Z_0 = X_0^2 - \sigma^2 + 2\lim_n X_0E_0(S_n), \qquad Z_l = Z_0\circ T^l \quad \text{and} \quad W_m = E_0(Z_1+Z_2+\cdots+Z_m).$$

*Suppose that $\sigma^2 > 0$. Let $\Lambda = \sigma^{-2}\mathbb{E}|X_0|^3$. Then*

$$d_1(S_n, \sigma\sqrt{n}Y) \le \frac{13\sigma}{6} + \frac{\Lambda}{6}\log(1+2n) + \sum_{m=1}^{[\sqrt{2n}]}\frac{\|X_0W_m\|_1 + 2\sigma\|W_m\|_1}{m\sigma^2} + D',$$

*where*

$$D' = \sum_{m=1}^{n}\frac{1}{\sigma\sqrt{m}}\left\|\sum_{l\ge m}X_0E_0(X_l)\right\|_1 + \sum_{m=1}^{n}\frac{1}{2m}\|(1+\sigma^{-2}X_0^2)E_0(S_m)\|_1.$$

***Remark 2.2.*** *By Theorem 1 in [8], the convergence in $\mathbb{L}^1$ of $X_0E_0(S_n)$ implies the convergence in distribution of $n^{-1/2}S_n$ to a mixture of Gaussian random variables. From [6] it also implies that $n^{-1/2}E_0(S_n)$ converges to 0 in $\mathbb{L}^1$. Consequently, if $X_0^2E_0(S_n)$ converges in $\mathbb{L}^1$ as $n\to\infty$, then $D' = \mathrm{o}(\sqrt{n})$. Moreover, from the $\mathbb{L}^1$-ergodic theorem, $(W_m/m)$ and $(X_0W_m/m)$ converge to 0 in $\mathbb{L}^1$ under the above additional condition. In that case, Theorem 2.2 gives a rate of convergence in the mean CLT.*

**Proof of Theorem 2.2.** Dividing the random variables by $\sigma$, we may assume that $\sigma = 1$. From (2.6) and and (2.9), for any 1-Lipschitzian function $f$,

$$\Delta(f) \le \sum_{k=1}^{n} \mathbb{E}\left( f_k'(S_{k-1})X_k + \frac{1}{2} f_k''(S_{k-1})(X_k^2 - 1) \right) + \rho(n), \tag{2.16}$$

where the functions $f_k$ are defined in Notation 2.1 of Section 2.1. In order to bound up the terms of first order, we write

$$f_k'(S_{k-1}) = f_0'(0) + \sum_{j=1}^{k-1} (f_{j+1}'(S_j) - f_j'(S_{j-1})).$$

Next

$$\begin{aligned} f_{j+1}'(S_j) - f_j'(S_{j-1}) &= (f_j'(S_j) - f_j'(S_{j-1})) - E_j(f_{j+1}'(S_j + Y) - f_{j+1}'(S_j)) \\ &= f_j''(S_{j-1})X_j + R_j', \end{aligned} \tag{2.17}$$

where $R_j'$ is some $\mathcal{M}_j$-measurable random variable such that

$$|R_j'| \le (2n - 2j)^{-1}(X_j^2 + 1). \tag{2.18}$$

Set $U_{j,n} = E_j(S_n - S_j)$. From (2.17) and (2.18)

$$\sum_{k=1}^{n} \mathbb{E}(f_k'(S_{k-1})X_k) \le \sum_{j=1}^{n-1} (\mathbb{E}(f_j''(S_{j-1})X_j U_{j,n}) + (2n - 2j)^{-1} \|(1 + X_j^2)U_{j,n}\|_1).$$

Next

$$\mathbb{E}(f_j''(S_{j-1})X_j U_{j,n}) \le \mathbb{E}(f_j''(S_{j-1})X_j U_{j,\infty}) + (n - j + 1)^{-1/2} \left\| \sum_{l>n} X_j E_j(X_l) \right\|_1$$

with the convention $X_j U_{j,\infty} = \lim_n X_j U_{j,n}$ in $\mathbb{L}^1$. From the stationarity, (2.16) and the above inequalities we get that

$$\Delta(f) \le \rho(n) + \frac{1}{2} \sum_{j=1}^{n} \mathbb{E}(f_j''(S_{j-1})Z_j) + D_1' + D_2', \tag{2.19}$$

where

$$D_1' = \sum_{m=1}^{n} m^{-1/2} \left\| \sum_{l \ge m} X_0 E_0(X_l) \right\|_1 \quad \text{and} \quad D_2' = \sum_{m=1}^{n} \frac{1}{2m} \|(1 + X_0^2)E_0(S_m)\|_1.$$

Theorem 2.2 follows then from (2.19) and Lemma 2.2 applied with $s = 2$. □

### *2.3. Applications to bounded random variables*

Throughout this subsection we assume that $X_0$ belongs to $\mathbb{L}^\infty$, and that $E_0(S_n)$ converges in $\mathbb{L}^1$. Then the series $X_0 E_0(S_n)$ converges in $\mathbb{L}^1$ and consequently Theorem 2.2 applies. Set

$$J_0 = \lim_{n\to\infty} E_0(S_n) \quad \text{and} \quad J_m = J_0 \circ T^m. \tag{2.20}$$

We first provide a rate which involves the quantities $\|E_0(m^{-1}S_m^2) - \sigma^2\|_1$ appearing in the conditional CLT of [6].

**Theorem 2.3.** *Let $(X_i)_{i\in\mathbb{Z}}$ be a stationary sequence of centered and bounded random variables such that $\mathbb{E}(X_0X_k|\mathcal{I})=\mathbb{E}(X_0X_k)$ a.s. for any integer $k$. Suppose furthermore that the sequence $E_0(S_n)$ converges in $\mathbb{L}^1$ to $J_0$. Then the series $\mathbb{E}(X_0^2)+2\sum_{k=1}^{\infty}\mathbb{E}(X_0X_k)$ is convergent to some nonnegative real $\sigma^2$ and $n^{-1}\operatorname{Var}S_n$ converges to $\sigma^2$. Suppose that $\sigma^2>0$ and let $L=\sigma^{-1}\|X_0\|_\infty$.*

(a) *If $S=\sum_{m\geq 0}\|E_0(J_m)\|_1<\infty$, then*

$$d_1(S_n,\sigma\sqrt{n}Y)\leq C\log(1+2n)+\sum_{m=1}^{[\sqrt{2n}]}\left(\frac{2+L}{m\sigma}\right)\|E_0(S_m^2)-m\sigma^2\|_1$$

*for some constant $C$ depending only on $\|X_0\|_\infty$, $\sigma$ and $S$.*

(b) *If $\|E_0(J_m)\|_1\leq M_\delta m^{-\delta}$ for some $\delta\in]0,1[$ and some constant $M_\delta$, then*

$$d_1(S_n,\sigma\sqrt{n}Y)\leq C_\delta n^{(1-\delta)/2}+\sum_{m=1}^{[\sqrt{2n}]}\left(\frac{2+L}{m\sigma}\right)\|E_0(S_m^2)-m\sigma^2\|_1$$

*for some constant $C_\delta$ depending on $\delta$, $M_\delta$, $\|X_0\|_\infty$ and $\sigma$.*

***Remark 2.3.*** *The assumptions made in this section ensure that*

$$\lim_{n\to\infty}\|E_0(m^{-1}S_m^2)-\sigma^2\|_1=0,$$

*which is the condition appearing in the conditional CLT of [6]. Consequently Theorem 2.3 provides rates of convergence in the mean CLT. For example, if* (a) *holds and $\sup_{m>0}\|E_0(S_m^2)-m\sigma^2\|_1<\infty$, then $d_1(S_n,\sigma\sqrt{n}Y)=\mathrm{O}(\log n)$. If* (b) *holds and $\|E_0(S_m^2)-m\sigma^2\|_1=\mathrm{O}(m^{1-\delta})$ as $m\to\infty$, then $d_1(S_n,\sigma\sqrt{n}Y)=\mathrm{O}(n^{(1-\delta)/2})$.*

**Proof of Theorem 2.3.** We first bound up $D'$. Let $M=\sup_{m>0}\|E_0(S_m)\|_1$. We have that

$$D'\leq L\sum_{m=1}^{n}m^{-1/2}\left\|\sum_{l\geq m}E_0(X_l)\right\|_1+\frac{1}{2}(1+L^2)M\log(1+2n).$$

Since $\sum_{l\geq m}E_0(X_l)=E_0(\sum_{l\geq m}E_{m-1}(X_l))=E_0(J_{m-1})$, we infer that

$$D'\leq L\sum_{m=0}^{n-1}(m+1)^{-1/2}\|E_0(J_m)\|_1+\frac{1}{2}(1+L^2)M\log(1+2n). \tag{2.21}$$

Next we bound up the r.v.'s $W_m+m\sigma^2-E_0(S_m^2)$ in $\mathbb{L}^1$. By definition of $W_m$,

$$W_m+m\sigma^2=E_0(S_m^2)+2\sum_{l=1}^{m}E_0\left(X_l\sum_{k>m}E_l(X_k)\right).$$

Therefore

$$\|W_m+m\sigma^2-E_0(S_m^2)\|_1\leq 2\sum_{l=1}^{m}\|X_lE_l(J_m)\|_1\leq 2\|X_0\|_\infty\sum_{l=1}^{m}\|E_0(J_{m-l})\|_1.$$

Hence

$$\sum_{m=1}^{[\sqrt{2n}]}(m\sigma^2)^{-1}(\|X_0W_m\|_1+\sigma\|W_m\|_1)$$

$$\leq \sum_{m=1}^{[\sqrt{2n}]} \left( \frac{2+L}{m\sigma} \right) \left( \|E_0(S_m^2) - m\sigma^2\|_1 + 2\|X_0\|_\infty \sum_{l=0}^{m-1} \|E_0(J_l)\|_1 \right). \tag{2.22}$$

Theorem 2.3 follows then easily from both Theorem 2.2, (2.21) and (2.22). □

We now give an application of Theorem 2.3 to sequences satisfying projective criteria in the style of [13, 14]. The proof, being elementary, is omitted.

**Corollary 2.1.** *Let $(X_i)_{i\in\mathbb{Z}}$ be a stationary sequence of centered and bounded random variables.*

(a) *If $\sum_{m=0}^{\infty}\sum_{l=0}^{m} \|E_{-m}(X_0X_l) - \mathbb{E}(X_0X_l)\|_1 < \infty$ and $\sum_{m>0} m\|E_0(X_m)\|_1 < \infty$, then the series of covariances converges to $\sigma^2$ and $d_1(S_n, \sigma\sqrt{n}Y) = \mathrm{O}(\log n)$ as $n$ tends to $\infty$, provided that $\sigma \neq 0$.*

(b) *If, for some $\delta \in ]0,1[$, $\sup_{l\in[0,m]} \|E_{-m}(X_0X_l) - \mathbb{E}(X_0X_l)\|_1 = \mathrm{O}(m^{-1-\delta})$ and $\|E_0(X_m)\|_1 = \mathrm{O}(m^{-1-\delta})$, then the series of covariances converges to $\sigma^2$ and $d_1(S_n, \sigma\sqrt{n}Y) = \mathrm{O}(n^{(1-\delta)/2})$ as $n$ tends to $\infty$ provided that $\sigma \neq 0$.*

***Remark 2.4.*** *For example, if the strong mixing coefficients $\alpha_2(k)$ of the sequence $(X_i)_{i\in\mathbb{Z}}$ (see (3.1) for the definition) satisfy $\alpha_2(k) = \mathrm{O}(k^{-1-\delta})$ then Corollary 2.1*(b) *applies and provides the rate of convergence $\mathrm{O}(n^{-\delta/2})$ in the mean CLT.*

## 3. Optimal rates for stationary sequences

Throughout Section 3, the filtration $(\mathcal{M}_i)_{i\in\mathbb{Z}}$ and the stationary sequence $(X_i)_{i\in\mathbb{Z}}$ are defined exactly as in Section 2.

### *3.1. Stationary sequences*

For stationary sequences, we will give two different conditions under which the rate of convergence $\mathrm{O}(n^{-1/2})$ holds in the mean CLT. We consider two types of dependence coefficients.

***Definition 3.1.*** *For any integers $0 \leq i < j$ and $p \geq 0$, let $\Gamma_{i,j,p}$ be the set of multiintegers $(k_1, \ldots, k_j)$ such that $0 \leq k_1 \leq \cdots \leq k_i$ and $k_i + p \leq k_{i+1} \leq \cdots \leq k_j$. Set*

$$\theta_{i,j}(p) = \sup_{(k_1,\ldots,k_j)\in\Gamma_{i,j,p}} \|X_{k_1}\cdots X_{k_i} E_{k_i}(X_{k_{i+1}}\cdots X_{k_j} - \mathbb{E}(X_{k_{i+1}}\cdots X_{k_j}))\|_1.$$

***Definition 3.2.*** *For any random variable $(\xi_1, \ldots, \xi_k)$ with values in $\mathbb{R}^k$, and any $\sigma$-algebra $\mathcal{M}$, define the function $g_{x,j}(t) = \mathbb{1}_{t\leq x} - \mathbb{P}(\xi_j \leq x)$. Set*

$$\alpha(\mathcal{M}, (\xi_1, \ldots, \xi_k)) = \sup_{(x_1,\ldots,x_k)\in\mathbb{R}^k} \left\| \mathbb{E}\left( \prod_{j=1}^{k} g_{x_j,j}(\xi_j) \Big| \mathcal{M} \right) - \mathbb{E}\left( \prod_{j=1}^{k} g_{x_j,j}(\xi_j) \right) \right\|_1.$$

*For a sequence $\boldsymbol{\xi} = (\xi_i)_{i\in\mathbb{Z}}$, where $\xi_i = \xi_0 \circ T^i$ and $\xi_0$ is a $\mathcal{M}_0$-measurable and real-valued r.v., let*

$$\alpha_{k,\boldsymbol{\xi}}(n) = \max_{1\leq l\leq k} \sup_{i_l>\cdots>i_1\geq n} \alpha(\mathcal{M}_0, (\xi_{i_1}, \ldots, \xi_{i_l})).$$

***Remark 3.1.*** *Let $B_1(\mathbb{R}^k)$ be the set of functions $f$ from $\mathbb{R}^k$ to $\mathbb{R}$ such that $|f(x) - f(y)| \leq 1$ for any $x, y$ in $\mathbb{R}^k$. Recall that the strong mixing coefficient of Rosenblatt may be defined as*

$$\alpha(\mathcal{M}, \sigma(\xi_1, \ldots, \xi_k)) = \frac{1}{2} \sup_{f\in B_1(\mathbb{R}^k)} \|\mathbb{E}(f(\xi_1, \ldots, \xi_k)|\mathcal{M}) - \mathbb{E}(f(\xi_1, \ldots, \xi_k))\|_1.$$

*For the sequence $\xi$, we define the strong mixing coefficients*

$$\alpha_k(n) = \sup_{i_k \geq \cdots \geq i_1 \geq n} \alpha(\mathcal{M}_0, \sigma(\xi_{i_1}, \dots, \xi_{i_l})) \quad \text{and} \quad \alpha(n) = \sup_{k>0} \alpha_k(n). \tag{3.1}$$

*By induction on $k$, it is easy to prove that $g:(t_1,\dots,t_k) \to \prod_{i=1}^k g_{x_i,i}(t_i)$ belongs to $B_1(\mathbb{R}^k)$. It follows that*

$$\alpha(\mathcal{M}, (\xi_1, \dots, \xi_k)) \leq 2\alpha(\mathcal{M}, \sigma(\xi_1, \dots, \xi_k)) \quad \text{and} \quad \alpha_{k,\boldsymbol{\xi}}(n) \leq 2\alpha_k(n).$$

*We emphasize that there exist sequences which are not strongly mixing in the sense of Rosenblatt, for which $\alpha_{k,\boldsymbol{\xi}}(n)$ tends to 0 as $n$ tends to infinity (see [7], Section 4 and the example of Section 5.1).*

***Definition 3.3.*** *For any real-valued random variable $X$, let $Q_X$ be the generalized inverse of the tail function $x \to \mathbb{P}(X > x)$.*

**Theorem 3.1.** *Let $(X_i)_{i\in\mathbb{Z}}$ be a stationary sequence of centered random variables. Consider the two conditions*

(a) $\mathbb{E}(X_0^4) < \infty$ *and* $\sum_{j=1}^{\infty} j\theta_{p,q}(j) < \infty$ *for any* $0 \leq p < q \leq (3+p) \wedge 4$.
(b) $X_0 = (f_1 - f_2)(\xi_0)$ *for some real-valued random variable $\xi_0$ and nondecreasing functions $f_1, f_2$, such that $f_1(\xi_0)$, $f_2(\xi_0)$ belong to $\mathbb{L}^3$ and, for $Q = \max(Q_{|f_1(\xi_0)|}, Q_{|f_2(\xi_0)|})$,*

$$\sum_{j=1}^{\infty} j \int_0^{\alpha_{3,\boldsymbol{\xi}}(j)} Q^3(u)\,\mathrm{d}u < \infty. \tag{3.2}$$

*If either* (a) *or* (b) *holds, then the series $\sigma^2 = \mathbb{E}(X_0^2) + 2\sum_{k=1}^{\infty} \mathbb{E}(X_0 X_k)$ converges absolutely. Moreover, if $\sigma > 0$, then $d_1(S_n, \sqrt{n}\sigma Y) \leq C$ for some constant $C$.*

***Remark 3.2.*** *For bounded random variables, Theorem 3.1 under* (a) *is a consequence of Theorem 1.1 in [22].*

***Remark 3.3.*** *For the strong mixing coefficients defined in (3.1), we infer from Theorem 3.1*(b) *that, if $X_0 = f(\xi_0)$ belongs to $\mathbb{L}^3$ and if (1.8) holds with $\alpha_3(k)$ instead of $\alpha(k)$ and $b = 1$, then the conclusion of Theorem 3.1 holds.*

### *3.2. Martingale difference sequences*

In this section we give conditions for stationary martingale difference sequences ensuring the optimal rate $\mathrm{O}(n^{-1/2})$ in the mean CLT.

**Theorem 3.2.** *Let $(X_i)_{i\in\mathbb{Z}}$ be a stationary martingale difference sequence in $\mathbb{L}^3$, with variance $\sigma^2$. Consider the two conditions*

(a) *$X_0$ belongs to $\mathbb{L}^4$,*

$$\sum_{k>0} \left( \|(X_0^2 \vee 1)(E_0(X_k^2) - \sigma^2)\|_1 + \frac{1}{k} \sum_{i=1}^{k} \|X_{-i} X_0 (E_0(X_k^2) - \sigma^2)\|_1 \right) < \infty,$$

*and*

$$\sum_{k>0} \frac{1}{k} \sum_{i=[k/2]}^{k} \|(|X_0| \vee 1)(E_0(X_i X_k^2) - \mathbb{E}(X_i X_k^2))\|_1 < \infty. \tag{3.3}$$

(b) *$X_0$ and $Q$ are defined as in Theorem 3.1*(b), *and*

$$\sum_{j=1}^{\infty} \int_0^{\alpha_{3,\boldsymbol{\xi}}(j)} Q^3(u)\,\mathrm{d}u < \infty. \tag{3.4}$$

*If either* (a) *or* (b) *holds, then $d_1(S_n, \sqrt{n}\sigma Y) \le C$ for some positive constant $C$.*

***Remark 3.4.*** *Note that the first condition in (3.3) implies that $\mathbb{E}(X_0^2|\mathcal{I}) = \sigma^2$ almost surely. Assume that $\mathbb{E}(|X_0|^p) < \infty$ for some $p \ge 4$. Applying Hölder's inequality, we see that (3.3) holds as soon as*

$$\sum_{k>0} \|E_0(X_k^2) - \sigma^2\|_{p/(p-2)} < \infty \quad \text{and} \quad \sum_{k>0} \sup_{i\ge k} \|E_0(X_k X_i^2) - \mathbb{E}(X_k X_i^2)\|_{p/(p-1)} < \infty.$$

## 4. Proofs of Theorems 3.1 and 3.2

### *4.1. A first decomposition*

The following proposition is the main step to prove Theorems 3.2 and 3.1. It is stated in the nonstationary case.

**Proposition 4.1.** *Let $(X_i)_{i\ge 1}$ be a sequence of centered random variables, each having a finite third moment, adapted to the filtration $(\mathcal{M}_i)_i$. Let $Z$ be a centered random variable with finite fourth moment independent of $\mathcal{M}_\infty$, and let $\mathbb{E}(Z^2) = \beta_2$, $\mathbb{E}(Z^3) = \beta_3$, $\mathbb{E}(Z^4) = \beta_4$. Let $S_0 = 0$ and $S_n = X_1 + \cdots + X_n$. Let $X_{i,1}$ and $X_{i,2}$ be two $\mathcal{M}_i$-measurable random variables such that $X_i = X_{i,1} + X_{i,2}$. For any four times continuously differentiable function $f$ and any $l$ in $[1, k[$,*

$$\mathbb{E}(f(S_{k-1} + X_k) - f(S_{k-1} + Z)) \le \zeta_1 A_1 + \zeta_2(A_2 + A_8) + \zeta_3(A_3 + A_9) + \zeta_4(A_4 + \cdots + A_7),$$

*where the reals $\zeta_i = \zeta_i(f)$ are defined by $\zeta_i = \|f^{(i)}\|_\infty$ and the numbers $A_i = A_i(k,l)$ are defined by*

$$A_1 = \|E_{k-l-1}(X_k)\|_1, \qquad A_2 = \frac{1}{2}\left\|\beta_2 - E_{k-l-1}\left(X_{k,1}X_k + 2\sum_{j=1}^{l} X_{k-j,1}X_k\right)\right\|_1,$$

$$A_3 = \frac{1}{6}\left\|\beta_3 - E_{k-l-1}\left(X_k X_{k,1}^2 + 3\sum_{j=1}^{l}(X_{k-j,1}(X_{k,1}X_k - \beta_2) + X_{k-j,1}^2 X_k)\right)\right.$$
$$\left. - 6E_{k-l-1}\left(\sum_{j=1}^{l}\sum_{p=1}^{j-1} X_{k-j,1}X_{k-p,1}X_k\right)\right\|_1,$$

$$A_4 = \frac{1}{24}(\mathbb{E}(|X_k X_{k,1}^3|) + \beta_4), \qquad A_5 = \frac{1}{6}\sum_{j=1}^{l}\|X_{k-j,1}^3 E_{k-j}(X_k)\|_1,$$

$$A_6 = \frac{1}{4}\sum_{j=1}^{l}\left\|X_{k-j,1}^2\left(\beta_2 - E_{k-j}\left(X_k X_{k,1} + 2\sum_{p=1}^{j-1} X_{k-p,1}X_k\right)\right)\right\|_1,$$

$$A_7 = \frac{1}{6}\sum_{j=1}^{l}\left\|X_{k-j}\left(\beta_3 - E_{k-j}\left(X_{k,1}^2 X_k + 3\sum_{p=1}^{j-1} X_{k-p,1}^2 X_k\right.\right.\right.$$
$$\left.\left.\left. + 3\sum_{p=1}^{j-1} X_{k-p,1}(X_{k,1}X_k - \beta_2) + 6\sum_{p=1}^{j-1}\sum_{q=1}^{p-1} X_{k-p,1}X_{k-q,1}X_k\right)\right)\right\|_1,$$

$$A_8 = \frac{1}{2}\left( \|X_{k,2}X_k\|_1 + 2\sum_{j=1}^{l} \left\| X_{k-j,2}E_{k-j}(X_k) \right\|_1 \right),$$

$$A_9 = \frac{1}{2}\sum_{j=1}^{l} \left\| X_{k-j,2}\left( \beta_2 - E_{k-j}\left( X_k X_{k,1} + 2\sum_{p=1}^{j-1} X_{k-p,1}X_k \right)\right)\right\|_1 .$$

**Proof.** We start from the equality

$$f(S_{k-1}+X_k) - f(S_{k-1}) = X_k \int_0^1 (f'(S_{k-1}+tX_{k,1}) - f'(S_{k-1}))\,\mathrm{d}t + X_k f'(S_{k-1}) + r_1(k),$$

with $r_1(k) \le (\zeta_2/2)|X_k X_{k,2}|$. Consequently

$$\begin{aligned} f(S_{k-1}+X_k) &= f(S_{k-1}) + f'(S_{k-1})X_k + \frac{X_k X_{k,1}}{2} f''(S_{k-1}) \\ &\quad + \frac{X_k X_{k,1}^2}{6} f'''(S_{k-1}) + R_1(k) + r_1(k), \end{aligned}$$

with $|R_1(k)| \le (\zeta_4/4!)|X_k X_{k,1}^3|$. Hence

$$\begin{aligned} \mathbb{E}(f(S_{k-1}+X_k) - f(S_{k-1}+Z)) &= \mathbb{E}(f'(S_{k-1})X_k) + \frac{1}{2}\mathbb{E}(f''(S_{k-1})(X_k X_{k,1} - \beta_2)) \\ &\quad + \frac{1}{6}\mathbb{E}(f'''(S_{k-1})(X_k X_{k,1}^2 - \beta_3)) + R_2(k) + \mathbb{E}(r_1(k)), \end{aligned}$$

with $R_2(k) \le \zeta_4 A_4$. Consider first the third-order terms. Clearly

$$\begin{aligned} \frac{1}{6} f'''(S_{k-1})(X_k X_{k,1}^2 - \beta_3) &= \frac{1}{6} f'''(S_{k-l-1})(X_k X_{k,1}^2 - \beta_3) \\ &\quad + \frac{1}{6}\sum_{j=1}^{l}\left(\int_0^1 f^{(4)}(S_{k-j-1}+tX_{k-j})\,\mathrm{d}t\right) X_{k-j}(X_k X_{k,1}^2 - \beta_3). \end{aligned} \tag{4.1}$$

Let $g_1(S_{k-j}, X_{k-j,1}) = f''(S_{k-j}) - f''(S_{k-j-1}+X_{k-j,1})$. For the second-order terms, we have first

$$\begin{aligned} \frac{1}{2} f''(S_{k-1})(X_k X_{k,1} - \beta_2) &= \frac{1}{2}(X_k X_{k,1} - \beta_2)\left( f''(S_{k-l-1}) + \sum_{j=1}^{l} g_1(S_{k-j}, X_{k-j,1}) + \sum_{j=1}^{l} f'''(S_{k-j-1})X_{k-j,1} \right. \\ &\quad \left. + \sum_{j=1}^{l}\left(\int_0^1 (1-t) f^{(4)}(S_{k-j-1}+tX_{k-j,1})\,\mathrm{d}t\right) X_{k-j,1}^2 \right), \end{aligned}$$

and next

$$\begin{aligned} \frac{1}{2} f''(S_{k-1})(X_k X_{k,1} - \beta_2) &= \frac{1}{2}(X_k X_{k,1} - \beta_2)\left( f''(S_{k-l-1}) + \sum_{j=1}^{l} g_1(S_{k-j}, X_{k-j,1}) + \sum_{j=1}^{l} f'''(S_{k-l-1})X_{k-j,1} \right. \\ &\quad + \sum_{j=1}^{l}\sum_{p=j+1}^{l}\left(\int_0^1 f^{(4)}(S_{k-p-1}+tX_{k-p})\,\mathrm{d}t\right) X_{k-p}X_{k-j,1} \\ &\quad \left. + \sum_{j=1}^{l}\left(\int_0^1 (1-t) f^{(4)}(S_{k-j-1}+tX_{k-j,1})\,\mathrm{d}t\right) X_{k-j,1}^2 \right). \end{aligned} \tag{4.2}$$

Let $g_2(S_{k-j}, X_{k-j,1}) = f'(S_{k-j}) - f'(S_{k-j-1} + X_{k-j,1})$. For the first-order terms, we have first

$$f'(S_{k-1})X_k = f'(S_{k-l-1})X_k + \sum_{j=1}^{l} g_2(S_{k-j}, X_{k-j,1})X_k + \sum_{j=1}^{l}(f'(S_{k-j-1} + X_{k-j,1}) - f'(S_{k-j-1}))X_k,$$

so that

$$\begin{aligned}f'(S_{k-1})X_k = f'(S_{k-l-1})X_k &+ \sum_{j=1}^{l} g_2(S_{k-j}, X_{k-j,1})X_k \\ &+ \sum_{j=1}^{l} f''(S_{k-j-1})X_{k-j,1}X_k + \frac{1}{2}\sum_{j=1}^{l} f'''(S_{k-j-1})X_{k-j,1}^2 X_k \\ &+ \sum_{j=1}^{l}\left(\int_0^1 \frac{(1-t)^2}{2} f^{(4)}(S_{k-j-1} + tX_{k-j,1})\,\mathrm{d}t\right)X_{k-j,1}^3 X_k.\end{aligned}$$

Next

$$\begin{aligned}f'(S_{k-1})X_k = X_k\Bigg(f'(S_{k-l-1}) &+ \sum_{j=1}^{l} g_2(S_{k-j}, X_{k-j,1}) + \sum_{j=1}^{l} f''(S_{k-l-1})X_{k-j,1} \\ &+ \frac{1}{2}\sum_{j=1}^{l} f'''(S_{k-j-1})X_{k-j,1}^2 + \sum_{j=1}^{l}\sum_{p=j+1}^{l} f'''(S_{k-p-1})X_{k-p,1}X_{k-j,1} \\ &+ \sum_{j=1}^{l}\sum_{p=j+1}^{l} X_{k-j,1}g_1(S_{k-p}, X_{k-p,1}) \\ &+ \sum_{j=1}^{l}\sum_{p=j+1}^{l}\left(\int_0^1 (1-t)f^{(4)}(S_{k-p-1} + tX_{k-p,1})\,\mathrm{d}t\right)X_{k-p,1}^2 X_{k-j,1} \\ &+ \sum_{j=1}^{l}\left(\int_0^1 \frac{(1-t)^2}{2} f^{(4)}(S_{k-j-1} + tX_{k-j,1})\,\mathrm{d}t\right)X_{k-j,1}^3\Bigg),\end{aligned}$$

whence

$$\begin{aligned}f'(S_{k-1})X_k = X_k\Bigg(f'(S_{k-l-1}) &+ \sum_{j=1}^{l} g_2(S_{k-j}, X_{k-j,1}) + \sum_{j=1}^{l} f''(S_{k-l-1})X_{k-j,1} \\ &+ \frac{1}{2}\sum_{j=1}^{l} f'''(S_{k-l-1})X_{k-j,1}^2 + \sum_{j=1}^{l}\sum_{p=j+1}^{l} f'''(S_{k-l-1})X_{k-p,1}X_{k-j,1} \\ &+ \sum_{j=1}^{l}\sum_{p=j+1}^{l} X_{k-j,1}g_1(S_{k-p}, X_{k-p,1}) \\ &+ \frac{1}{2}\sum_{j=1}^{l}\sum_{p=j+1}^{l}\left(\int_0^1 f^{(4)}(S_{k-p-1} + tX_{k-p})\,\mathrm{d}t\right)X_{k-p}X_{k-j,1}^2\end{aligned}$$

$$+\sum_{j=1}^{l}\sum_{p=j+1}^{l}\sum_{q=p+1}^{l}\left(\int_0^1 f^{(4)}(S_{k-q-1}+tX_{k-q})\,\mathrm{d}t\right)X_{k-q}X_{k-p,1}X_{k-j,1}$$

$$+\sum_{j=1}^{l}\sum_{p=j+1}^{l}\left(\int_0^1 (1-t)f^{(4)}(S_{k-p-1}+tX_{k-p,1})\,\mathrm{d}t\right)X_{k-p,1}^2X_{k-j,1}$$

$$+\sum_{j=1}^{l}\left(\int_0^1 \frac{(1-t)^2}{2}f^{(4)}(S_{k-j-1}+tX_{k-j,1})\,\mathrm{d}t\right)X_{k-j,1}^3\Bigg). \tag{4.3}$$

Let us look carefully at the decompositions (4.1), (4.2) and (4.3). In front of $f'(S_{k-l-1})$ there is $X_k$, which leads to the term $\zeta_1 A_1$ by taking the conditional expectation with respect to $\mathcal{M}_{k-l-1}$. In front of $f''(S_{k-l-1})/2$ there is $X_kX_{k,1}-\beta_2+2\sum_{j=1}^{l}X_{k-j,1}X_k$, which leads to the term $\zeta_2A_2$ by taking the conditional expectation with respect to $\mathcal{M}_{k-l-1}$. In front of $f'''(S_{k-l-1})/6$ there is

$$X_kX_{k,1}^2-\beta_3+3\sum_{j=1}^{l}(X_{k-j,1}(X_kX_{k,1}-\beta_2)+X_{k-j,1}^2X_k)+6\sum_{p=1}^{l}\sum_{j=1}^{p-1}X_{k-p,1}X_{k-j,1}X_k,$$

which leads to the term $\zeta_3A_3$ by taking the conditional expectation with respect to $\mathcal{M}_{k-l-1}$. Taking the conditional expectation with respect to $\mathcal{M}_{k-j}$ and the supremum of $|f^{(4)}|$ in the last term of (4.3), we obtain $\zeta_4A_5$. Gathering the last term in (4.2) and the last but one in (4.3), we obtain

$$\sum_{j=1}^{l}\left(\int_0^1(1-t)f^{(4)}(S_{k-j-1}+tX_{k-j,1})\,\mathrm{d}t\right)X_{k-j,1}^2\left(\frac{1}{2}(X_{k,1}X_k-\beta_2)+\sum_{p=1}^{j-1}X_{k-p,1}X_k\right),$$

which leads to the term $\zeta_4A_6$. Gathering the remainder terms in (4.1), (4.2) and (4.3) (except the terms involving the functions $g_1,g_2$), we obtain

$$\sum_{j=1}^{l}\left(\int_0^1 f^{(4)}(S_{k-j-1}+tX_{k-j})\,\mathrm{d}t\right)X_{k-j}\Bigg(\frac{1}{6}(X_kX_{k,1}^2-\beta_3)+\frac{1}{2}\sum_{p=1}^{j-1}X_{k-p,1}^2X_k$$
$$+\frac{1}{2}\sum_{p=1}^{j-1}X_{k-p,1}(X_{k,1}X_k-\beta_2)+\sum_{p=1}^{j-1}\sum_{q=1}^{p-1}X_{k-p,1}X_{k-q,1}X_k\Bigg),$$

which leads to the term $\zeta_4A_7$. The term $\zeta_2A_8$ is obtained by gathering $\|r_1(k)\|_1$ and the terms involving the function $g_2$, and by noting that $|g_2(S_{k-j},X_{k-j,1})|\le\zeta_2|X_{k-j,2}|$. The term $\zeta_3A_9$ is obtained by gathering the terms involving the function $g_1$, and by noting that $|g_1(S_{k-j},X_{k-j,1})|\le\zeta_3|X_{k-j,2}|$. □

### *4.2. Upper bounds for the $A_i$'s*

Let $X_{i,1}$ and $X_{i,2}$ be two $\mathcal{M}_i$-measurable random variables such that $X_i=X_{i,1}+X_{i,2}$. Define $b(l)$ by

$$b(l)=\mathbb{E}(X_{0,1}^2X_0)+3\sum_{i=1}^{l}\mathbb{E}(X_{0,1}X_{i,1}X_i+X_{0,1}^2X_i)+6\sum_{i=1}^{l}\sum_{j=1}^{i-1}\mathbb{E}(X_{0,1}X_{j,1}X_i). \tag{4.4}$$

Assume that the series

$$\sigma^2=\mathbb{E}(X_0^2)+2\sum_{k=1}^{\infty}\mathbb{E}(X_0X_k)\quad\text{and}\quad\sigma_1^2=\mathbb{E}(X_{0,1}X_0)+2\sum_{k=1}^{\infty}\mathbb{E}(X_{0,1}X_k)$$

converge absolutely. Let $A_i$ be the terms of Proposition 4.1 with $\beta_2 = \sigma^2$, and $\beta_3 = b(l)$, and let $A_{i,1}$ be the terms of Proposition 4.1 with $\beta_2 = \sigma_1^2$, and $\beta_3 = b(l)$. We now give upper bounds for $A_{2,1}, A_{3,1}, A_{6,1}, A_{7,1}$ and $A_{9,1}$. First,

$$A_{2,1} \le \frac{1}{2}\|E_{k-l-1}(X_{k,1}X_k - \mathbb{E}(X_{k,1}X_k))\|_1 + \sum_{j=[l/2]+1}^{l} \|X_{k-j,1}E_{k-j}(X_k)\|_1 + \sum_{j=[l/2]+1}^{\infty} |\mathbb{E}(X_{0,1}X_j)| + \sum_{j=1}^{[l/2]} \|E_{k-l-1}(X_{k-j,1}X_k - \mathbb{E}(X_{k-j,1}X_k))\|_1,$$

$$2A_{6,1} \le \sum_{j=1}^{l}\Bigg(\frac{1}{2}\|X_{k-j,1}^2 E_{k-j}(X_{k,1}X_k - \mathbb{E}(X_{k,1}X_k))\|_1 + \sum_{p=[j/2]+1}^{j} \|X_{k-j,1}^2 X_{k-p,1}E_{k-p}(X_k)\|_1 + \mathbb{E}(X_{0,1}^2)\sum_{p=[j/2]+1}^{\infty} |\mathbb{E}(X_{0,1}X_p)| + \sum_{p=1}^{[j/2]} \|X_{k-j,1}^2 E_{k-j}(X_{k-p,1}X_k - \mathbb{E}(X_{k-p,1}X_k))\|_1\Bigg),$$

$$A_{9,1} \le \sum_{j=1}^{l}\Bigg(\frac{1}{2}\|X_{k-j,2} E_{k-j}(X_{k,1}X_k - \mathbb{E}(X_{k,1}X_k))\|_1 + \sum_{p=[j/2]+1}^{j} \|X_{k-j,2} X_{k-p,1}E_{k-p}(X_k)\|_1 + \|X_{0,2}\|_1\sum_{p=[j/2]+1}^{\infty} |\mathbb{E}(X_{0,1}X_p)| + \sum_{p=1}^{[j/2]} \|X_{k-j,2} E_{k-j}(X_{k-p,1}X_k - \mathbb{E}(X_{k-p,1}X_k))\|_1\Bigg).$$

Next, we have that $A_{3,1} \le C_1 + C_2 + C_3$, where

$$C_1 = \frac{1}{6}\Bigg(\|E_{k-l-1}(X_{k,1}^2X_k - \mathbb{E}(X_{k,1}^2X_k))\|_1 + 3\sum_{j=[l/2]+1}^{l} \|X_{k-j,1}^2E_{k-j}(X_k)\|_1 + 3\sum_{j=1}^{[l/2]} \|E_{k-l-1}(X_{k-j,1}^2X_k - \mathbb{E}(X_{k-j,1}^2X_k))\|_1 + 3\sum_{j=[l/2]+1}^{l} \|X_{k-j,1}E_{k-j}(X_{k,1}X_k - \mathbb{E}(X_{k,1}X_k))\|_1 + 3\sum_{j=1}^{[l/2]} \|E_{k-l-1}(X_{k-j,1}X_{k,1}X_k - \mathbb{E}(X_{k-j,1}X_{k,1}X_k))\|_1\Bigg),$$

$$C_2 = \sum_{j=1}^{[l/2]}\sum_{p=1}^{j-1} \|E_{k-l-1}(X_{k-j,1}X_{k-p,1}X_k - \mathbb{E}(X_{k-j,1}X_{k-p,1}X_k))\|_1 + \sum_{j=[l/2]+1}^{l}\sum_{p=1}^{[j/2]} \|X_{k-j,1}E_{k-j}(X_{k-p,1}X_k - \mathbb{E}(X_{k-p,1}X_k))\|_1 + \sum_{j=[l/2]+1}^{l}\sum_{p=[j/2]+1}^{j-1} \|X_{k-j,1}X_{k-p,1}E_{k-p}(X_k)\|_1,$$

$$C_3 = \sum_{j=[l/2]+1}^{l}\sum_{p=1}^{j-1} |\mathbb{E}(X_{0,1}X_{p,1}X_j)| + \frac{1}{2}\sum_{j=[l/2]+1}^{l} |\mathbb{E}(X_{0,1}^2X_j)|$$

$$
\begin{aligned}
&+\frac{1}{2}\sum_{j=1}^{[l/2]}\|E_{k-l-1}(X_{k-j,1})\|_1\biggl(|\mathbb{E}(X_{0,1}X_0)|+2\sum_{p=1}^{\infty}|\mathbb{E}(X_{0,1}X_p)|\biggr)\\
&+\frac{1}{2}\sum_{j=[l/2]+1}^{l}|\mathbb{E}(X_{0,1}X_{j,1}X_j)|+\sum_{j=[l/2]+1}^{l}\sum_{p=[j/2]+1}^{\infty}\|X_{0,1}\|_1|\mathbb{E}(X_{0,1}X_p)|.
\end{aligned}
$$

In the same way, $A_{7,1}\le D_1+D_2+D_3$, where

$$
\begin{aligned}
D_1=\frac{1}{6}\sum_{j=1}^{l}\biggl(&\|X_{k-j}E_{k-j}(X_{k,1}^2X_k-\mathbb{E}(X_{k,1}^2X_k))\|_1\\
&+3\sum_{p=[j/2]+1}^{j-1}\|X_{k-j}X_{k-p,1}^2E_{k-p}(X_k)\|_1\\
&+3\sum_{p=1}^{[j/2]}\|X_{k-j}E_{k-j}(X_{k-p,1}^2X_k-\mathbb{E}(X_{k-p,1}^2X_k))\|_1\\
&+3\sum_{p=[j/2]+1}^{j-1}\|X_{k-j}X_{k-p,1}E_{k-p}(X_{k,1}X_k-\mathbb{E}(X_{k,1}X_k))\|_1\\
&+3\sum_{p=1}^{[j/2]}\|X_{k-j}E_{k-j}(X_{k-p,1}X_{k,1}X_k-\mathbb{E}(X_{k-p,1}X_{k,1}X_k))\|_1\biggr),\\
D_2=\sum_{j=1}^{l}\biggl(&\sum_{p=[j/2]+1}^{j-1}\sum_{q=[p/2]+1}^{p-1}\|X_{k-j}X_{k-p,1}X_{k-q,1}E_{k-q}(X_k)\|_1\\
&+\sum_{p=1}^{[j/2]}\sum_{q=1}^{p-1}\|X_{k-j}E_{k-j}(X_{k-p,1}X_{k-q,1}X_k-\mathbb{E}(X_{k-p,1}X_{k-q,1}X_k))\|_1\\
&+\sum_{p=[j/2]+1}^{j-1}\sum_{q=1}^{[p/2]}\|X_{k-j}X_{k-p,1}E_{k-p}(X_{k-q,1}X_k-\mathbb{E}(X_{k-q,1}X_k))\|_1\biggr),\\
D_3=\sum_{j=1}^{l}\biggl(&2\sum_{p=[j/2]+1}^{l}\sum_{q=1}^{p-1}\|X_0\|_1|\mathbb{E}(X_{0,1}X_{q,1}X_p)|\\
&+\sum_{p=1}^{[j/2]}\|X_{k-j}E_{k-j}(X_{k-p,1})\|_1\biggl(|\mathbb{E}(X_0X_{0,1})|+2\sum_{p=1}^{\infty}|\mathbb{E}(X_{0,1}X_p)|\biggr)\\
&+\sum_{p=[j/2]+1}^{l}\|X_0\|_1|\mathbb{E}(X_{0,1}^2X_p)|+\sum_{p=[j/2]+1}^{l}\|X_0\|_1|\mathbb{E}(X_{0,1}X_{p,1}X_p)|\\
&+2\sum_{p=[j/2]+1}^{l}\sum_{q=[p/2]+1}^{\infty}\|X_0\|_2\|X_{0,1}\|_2|\mathbb{E}(X_{0,1}X_q)|\biggr).
\end{aligned}
$$

### *4.3. Control of the $A_i$'s for stationary sequences*

In this section, we give bounds for the quantities $A_i$ for $i \neq 4$. The control of $A_4$ is carried out in Section 4.5. The bounds are given in terms of the coefficients $\theta$, and in terms of $\alpha_\xi$ in the case where $X_0 = (f_1 - f_2)(\xi_0)$, the functions $f_1, f_2$ being nondecreasing. For $\alpha_\xi$, let

$$X_{i,1} = X_i(a) = (g_a \circ f_1 - g_a \circ f_2)(\xi_i) - \mathbb{E}((g_a \circ f_1 - g_a \circ f_2)(\xi_i)),$$

where $g_a(x) = (x \wedge a) \vee (-a)$ for any $a > 0$. For $\theta$, let $X_{i,1} = X_i(\infty) = X_i$, in which case $A_{i,1} = A_i$ and $A_8(f,k,l) = A_9(f,k,l) = 0$. Denote by $b(l,a)$ the quantity $b(l)$ defined in (4.4) with $X_{i,1} = X_i(a)$. Note that $b(l,\infty)$ converges to a limit $b(\infty,\infty)$ as soon as both $\sum k\theta_{1,3}(k)$ and $\sum k\theta_{2,3}(k)$ are finite. In the same way, since $g_a \circ f_1$ and $g_a \circ f_2$ are nondecreasing, we can use Corollary A.1 given in the Appendix: it follows easily that $b(l,a)$ converges to a limit $b(\infty,a)$ as soon as (3.2) holds.

***Notation 4.1.*** *In the following, the notation $a \ll b$ means that $a \le Cb$ for some numerical constant $C$.*

To control the $A_i$'s with the help of the coefficients $\alpha_{\boldsymbol{\xi}}$, the main tool is the second inequality given in Corollary A.1 of the Appendix. Let $X_k^{(a)} = X_{k,2} = X_k - X_k(a)$. Then $\sigma^2 - \sigma_1^2 = \mathbb{E}(X_0 X_0^{(a)}) + 2\sum_{k=1}^{\infty} \mathbb{E}(X_0^{(a)} X_k)$. Note that

$$X_k^{(a)} = (h_a \circ f_1 - h_a \circ f_2)(\xi_k) - \mathbb{E}((h_a \circ f_1 - h_a \circ f_2)(\xi_k)),$$

where $h_a(x) = x - g_a(x)$. The functions $h_a \circ f_1$ and $h_a \circ f_2$ are nondecreasing and

$$\max(Q_{|g_a \circ f_1(Y_0)|}, Q_{|g_a \circ f_2(Y_0)|}, Q_{|h_a \circ f_1(Y_0)|}, Q_{|h_a \circ f_2(Y_0)|}) \le Q.$$

Hence Corollary A.1 applies and yields $|\sigma^2 - \sigma_1^2| \ll M(Q,a)$ where

$$M(Q,a) = \sum_{i=0}^{\infty} \int_0^{\alpha_{1,\xi}(i)} Q^2 \mathbb{1}_{Q>a} \, d\lambda,$$

$\lambda$ being the Lebesgue measure. Taking into account this upper bound, we get that

$$\begin{aligned}
&|A_2 - A_{2,1}| \ll M(Q,a),\\
&|A_3 - A_{3,1}| \ll l\|X_0(a)\|_1 M(Q,a),\\
&|A_6 - A_{6,1}| \ll l\|X_0(a)\|_2^2 M(Q,a),\\
&|A_7 - A_{7,1}| \ll l^2\|X_0\|_2^2 M(Q,a),\\
&|A_9 - A_{9,1}| \ll l\|X_0^{(a)}\|_1 M(Q,a).
\end{aligned} \tag{4.5}$$

We now give some upper bounds for the $A_i$'s. Clearly

$$A_1 \le \theta_{0,1}(l+1) \quad \text{and} \quad A_1 \ll \int_0^{\alpha_{1,\xi}(l)} Q \, d\lambda. \tag{4.6}$$

In the same way, since $\max(Q_{|g_a \circ f_1(Y_0)|}, Q_{|g_a \circ f_2(Y_0)|}) \le (Q \wedge a)$,

$$A_5 \ll \sum_{j=1}^{l} \theta_{3,4}(j) \quad \text{and} \quad A_5 \ll \sum_{j=1}^{l} \int_0^{\alpha_{1,\xi}(j)} Q(Q \wedge a)^3 \, d\lambda. \tag{4.7}$$

Let $A = \operatorname{sign}\{E_{k-l-1}(X_k(a)X_k - \mathbb{E}(X_k(a)X_k)\}$. Recall that $\alpha(\xi_1,\dots,\xi_k)$ is defined in Proposition A.1. Since $\alpha(A,\xi_k,\xi_k) \le \alpha(A,\xi_k) \le \alpha_{1,\boldsymbol{\xi}}(l+1)$, we infer from Corollary A.1 that

$$\|E_{k-l-1}(X_k(a)X_k - \mathbb{E}(X_k(a)X_k))\|_1 = |\mathbb{E}((A-\mathbb{E}(A))X_k(a)X_k)| \ll \int_0^{\alpha_{1,\boldsymbol{\xi}}(l+1)} Q^2\,\mathrm{d}\lambda.$$

Using this inequality to control $A_{2,1}$, we obtain the bounds

$$A_2 \ll \theta_{0,2}(l+1) + \sum_{j=[l/2]}^{\infty} \theta_{0,2}(j) + \sum_{j=[l/2]}^{\infty} \theta_{1,2}(j), \tag{4.8}$$

$$A_{2,1} \ll \int_0^{\alpha_{1,\boldsymbol{\xi}}(l+1)} Q^2\,\mathrm{d}\lambda + \sum_{j=[l/2]}^{\infty} \int_0^{\alpha_{2,\boldsymbol{\xi}}(j)} Q^2\,\mathrm{d}\lambda. \tag{4.9}$$

In the same way,

$$A_6 \ll \sum_{j=1}^{\infty} j\theta_{2,4}(j) + \sum_{p=1}^{\infty} p\theta_{3,4}(p) + \mathbb{E}(X_0^2)\sum_{p=1}^{\infty} p\theta_{1,2}(p), \tag{4.10}$$

$$\begin{aligned} A_{6,1} \ll{}& \sum_{j=1}^{l} \int_0^{\alpha_{1,\boldsymbol{\xi}}(j)} Q(Q\wedge a)^3\,\mathrm{d}\lambda + \sum_{p=1}^{l} p\int_0^{\alpha_{2,\boldsymbol{\xi}}(p)} Q(Q\wedge a)^3\,\mathrm{d}\lambda \\ &+ \mathbb{E}(X_0^2(a))\sum_{p=1}^{\infty} p\int_0^{\alpha_{1,\boldsymbol{\xi}}(p)} Q(Q\wedge a)\,\mathrm{d}\lambda. \end{aligned} \tag{4.11}$$

The term $A_{9,1}$ can be handled similarly:

$$\begin{aligned} A_{9,1} \ll{}& a\sum_{j=1}^{l} \int_0^{\alpha_{1,\boldsymbol{\xi}}(j)} Q^2 \mathbb{1}_{Q>a}\,\mathrm{d}\lambda + a\sum_{p=1}^{l} p\int_0^{\alpha_{2,\boldsymbol{\xi}}(p)} Q^2\mathbb{1}_{Q>a}\,\mathrm{d}\lambda \\ &+ \|X_0^{(a)}\|_1 \sum_{p=1}^{\infty} p\int_0^{\alpha_{1,\boldsymbol{\xi}}(p)} Q(Q\wedge a)\,\mathrm{d}\lambda. \end{aligned} \tag{4.12}$$

In the previous section, we have defined quantities $C_1, C_2, C_3$ such that $A_{3,1} \le C_1 + C_2 + C_3$. If $X_{k,1} = X_k(a)$ we shall use the notation $C_i = C_i(a)$, and if $X_{k,1} = X_k(\infty) = X_k$ the notation $C_i = C_i(\infty)$. Thus $A_{3,1} \le C_1(a) + C_2(a) + C_3(a)$ and $A_3 \le C_1(\infty) + C_2(\infty) + C_3(\infty)$. To control $C_i(a)$, we use Corollary A.1 and the fact that, for any $\mathcal{M}_0$-measurable r.v. $B$,

$$\alpha(B,\xi_k,\xi_k,\xi_k) \le \alpha_{1,\boldsymbol{\xi}}(k) \quad \text{and} \quad \alpha(B,\xi_k,\xi_k,\xi_l) \le \alpha_{2,\boldsymbol{\xi}}(\min(k,l)).$$

Therefrom

$$C_1(\infty) \ll \sum_{j=[l/2]}^{l+1} \theta_{2,3}(j) + 2\sum_{j=[l/2]}^{l} \theta_{0,3}(j) + \sum_{j=[l/2]}^{l} \theta_{1,3}(j), \tag{4.13}$$

$$C_1(a) \ll \int_0^{\alpha_{1,\boldsymbol{\xi}}(l+1)} Q(Q\wedge a)^2\,\mathrm{d}\lambda + \sum_{j=[l/2]}^{l} \int_0^{\alpha_{2,\boldsymbol{\xi}}(j)} Q(Q\wedge a)^2\,\mathrm{d}\lambda, \tag{4.14}$$

$$C_2(\infty) \ll l\left(\sum_{j=[l/2]}^{l} \theta_{0,3}(j) + \sum_{j=[l/4]}^{l} \theta_{1,3}(j) + \sum_{j=[l/4]}^{l} \theta_{2,3}(j)\right), \tag{4.15}$$

$$C_2(a) \ll l \sum_{j=[l/4]}^{l} \int_0^{\alpha_{3,\boldsymbol{\xi}}(j)} Q(Q\wedge a)^2 \, \mathrm{d}\lambda. \tag{4.16}$$

Finally

$$\begin{aligned} C_3(\infty) \ll{}& l \sum_{j=[l/4]}^{l} (\theta_{1,3}(j)+\theta_{2,3}(j)) + \sum_{j=[l/2]}^{\infty} (\theta_{1,3}(j)+\theta_{2,3}(j)) \\ &+ \left(\sum_{j=[l/2]}^{l} \theta_{0,1}(j)\right)\left(\mathbb{E}(X_0^2) + 2\sum_{p=1}^{\infty} |\mathbb{E}(X_0X_p)|\right) + l\|X_0\|_1 \sum_{j=[l/4]}^{\infty} \theta_{1,2}(j) \end{aligned} \tag{4.17}$$

and

$$\begin{aligned} C_3(a) \ll{}& l \sum_{j=[l/4]}^{l} \int_0^{\alpha_{2,\boldsymbol{\xi}}(j)} Q(Q\wedge a)^2 \, \mathrm{d}\lambda + \sum_{j=[l/2]}^{\infty} \int_0^{\alpha_{1,\boldsymbol{\xi}}(j)} Q(Q\wedge a)^2 \, \mathrm{d}\lambda \\ &+ \left(\sum_{j=[l/2]}^{l} \int_0^{\alpha_{1,\boldsymbol{\xi}}(j)} Q \, \mathrm{d}\lambda\right)\left(|\mathbb{E}(X_0X_0(a))| + 2\sum_{p=1}^{\infty} |\mathbb{E}(X_0(a)X_p)|\right) \\ &+ l\|X_0(a)\|_1 \sum_{j=[l/4]}^{\infty} \int_0^{\alpha_{1,\boldsymbol{\xi}}(j)} Q(Q\wedge a) \, \mathrm{d}\lambda. \end{aligned} \tag{4.18}$$

In the previous section, we have defined quantities $D_1, D_2, D_3$ such that $A_{7,1} \le D_1 + D_2 + D_3$. If $X_{k,1} = X_k(a)$ we shall use the notation $D_i = D_i(a)$, and if $X_{k,1} = X_k(\infty) = X_k$ the notation $D_i = D_i(\infty)$. Thus $A_{7,1} \le D_1(a) + D_2(a) + D_3(a)$ and $A_7 \le D_1(\infty) + D_2(\infty) + D_3(\infty)$. Furthermore

$$D_1(\infty) \ll \sum_{j=1}^{l} \theta_{1,4}(j) + \sum_{j=1}^{l} j(\theta_{3,4}(j) + \theta_{1,4}(j) + \theta_{2,4}(j)), \tag{4.19}$$

$$D_1(a) \ll \sum_{j=1}^{l} \int_0^{\alpha_{1,\boldsymbol{\xi}}(j)} Q^2(Q\wedge a)^2 \, \mathrm{d}\lambda + \sum_{j=1}^{l} j \int_0^{\alpha_{2,\boldsymbol{\xi}}(j)} Q^2(Q\wedge a)^2 \, \mathrm{d}\lambda, \tag{4.20}$$

$$D_2(\infty) \ll \sum_{j=1}^{l} (l\wedge 2j)^2 \theta_{1,4}(j) + \sum_{j=1}^{l} (l\wedge 2j)^2 (\theta_{2,4}(j) + \theta_{3,4}(j)), \tag{4.21}$$

$$D_2(a) \ll \sum_{j=1}^{l} (l\wedge 2j)^2 \int_0^{\alpha_{3,\boldsymbol{\xi}}(j)} Q^2(Q\wedge a)^2 \, \mathrm{d}\lambda. \tag{4.22}$$

Finally

$$\begin{aligned} D_3(\infty) \ll{}& \|X_0\|_1 \sum_{j=1}^{l} (l\wedge 2j)^2 (\theta_{1,3}(j) + \theta_{2,3}(j)) + \|X_0\|_2^2 \sum_{j=1}^{\infty} (l\wedge 2j)^2 \theta_{1,2}(j) \\ &+ \left(\sum_{j=1}^{l} j\theta_{1,2}(j)\right)\left(\mathbb{E}(X_0^2) + 2\sum_{p=1}^{\infty} |\mathbb{E}(X_0X_p)|\right) \end{aligned} \tag{4.23}$$

and

$$D_3(a) \ll \|X_0\|_1 \left( \sum_{j=1}^{l} (l \wedge 2j)^2 \int_0^{\alpha_{2,\boldsymbol{\xi}}(j)} Q(Q \wedge a)^2 \, \mathrm{d}\lambda + \sum_{j=1}^{\infty} j \int_0^{\alpha_{1,\boldsymbol{\xi}}(j)} Q(Q \wedge a)^2 \, \mathrm{d}\lambda \right)$$
$$+ \left( \sum_{j=1}^{l} j \int_0^{\alpha_{1,\boldsymbol{\xi}}(j)} Q^2 \, \mathrm{d}\lambda \right) \left( |\mathbb{E}(X_0 X_0(a))| + 2 \sum_{p=1}^{\infty} |\mathbb{E}(X_0(a) X_p)| \right)$$
$$+ \|X_0\|_2 \|X_0(a)\|_2 \sum_{j=1}^{\infty} (l \wedge 2j)^2 \int_0^{\alpha_{1,\boldsymbol{\xi}}(j)} Q(Q \wedge a) \, \mathrm{d}\lambda. \tag{4.24}$$

It remains to bound up $A_8$. Clearly

$$A_8 \ll \sum_{j=0}^{l} \int_0^{\alpha_{1,\boldsymbol{\xi}}(j)} Q^2 \mathbb{1}_{Q>a} \, \mathrm{d}\lambda. \tag{4.25}$$

### *4.4. Control of the $A_i$'s for martingales*

For stationary martingale difference sequences, the control of the eight terms $A_i$ is much easier, since the terms $A_1, A_5, C_2, D_2$ are equal to 0. If moreover $X_{k,1} = X_k(\infty) = X_k$, then $A_8$ and $A_9$ are equal to 0. We start from the control of the previous section. For each term $A_i$, we shall first give an upper bound when $X_{k,1} = X_k(\infty) = X_k$ in terms of sums of conditional expectations, and next an upper bound involving the mixing coefficients $\alpha_{\boldsymbol{\xi}}$. Clearly,

$$A_2 \le \frac{1}{2} \|E_0(X_l^2) - \sigma^2\|_1 \quad \text{and} \quad A_{2,1} \ll \int_0^{\alpha_{1,\boldsymbol{\xi}}(l+1)} Q^2 \, \mathrm{d}\lambda.$$

In the same way

$$A_6 \le \frac{1}{4} \sum_{j=1}^{l} \|X_0^2 (E_0(X_j^2) - \sigma^2)\|_1, \quad A_{6,1} \ll \sum_{j=1}^{l} \int_0^{\alpha_{2,\boldsymbol{\xi}}(j)} Q(Q \wedge a)^3 \, \mathrm{d}\lambda$$

and

$$A_{9,1} \ll a \sum_{j=1}^{l} \int_0^{\alpha_{2,\boldsymbol{\xi}}(j)} Q^2 \mathbb{1}_{Q>a} \, \mathrm{d}\lambda.$$

Starting from the control $A_{3,1} \le C_1 + C_2 + C_3$, and noting that $C_2 = 0$ for martingale difference sequences, we infer that

$$A_3 \ll \|E_0(X_l^3) - \mathbb{E}(X_l^3)\|_1 + \sum_{j=[l/2]}^{l} \|X_0(E_0(X_j^2) - \sigma^2)\|_1 + \sum_{j=[l/2]}^{l-1} \|E_0(X_j X_l^2) - \mathbb{E}(X_j X_l^2)\|_1,$$

$$A_{3,1} \ll \int_0^{\alpha_{1,\boldsymbol{\xi}}(l+1)} Q(Q \wedge a)^2 \, \mathrm{d}\lambda + \sum_{j=[l/2]}^{l} \int_0^{\alpha_{2,\boldsymbol{\xi}}(j)} Q(Q \wedge a)^2 \, \mathrm{d}\lambda + |\mathbb{E}(X_0 X_0(a))| \sum_{j=[l/2]}^{l} \int_0^{\alpha_{1,\boldsymbol{\xi}}(j)} Q \, \mathrm{d}\lambda.$$

Starting from the control $A_{7,1} \le D_1(\infty) + D_2(\infty) + D_3(\infty)$, and noting that $D_2(\infty) = 0$ for martingale difference sequences, we infer that

$$A_7 \ll \sum_{j=1}^{l} \Bigl( \|X_0(E_0(X_j^3) - \mathbb{E}(X_j^3))\|_1 + j \|X_0\|_1 \|X_0(E_0(X_j^2) - \sigma^2)\|_1$$

$$+\sum_{p=1}^{j}\|X_{-p}X_0(E_0(X_j^2)-\sigma^2)\|_1+\sum_{p=[j/2]}^{j-1}\|X_0(E_0(X_pX_j^2)-\mathbb{E}(X_pX_j^2))\|_1\Bigg),$$

$$A_{7,1}\ll\sum_{j=1}^{l}\int_0^{\alpha_{1,\xi}(j)}Q^2(Q\wedge a)^2\,\mathrm{d}\lambda+\sum_{j=1}^{l}j\int_0^{\alpha_{2,\xi}(j)}Q^2(Q\wedge a)^2\,\mathrm{d}\lambda$$

$$+\|X_0\|_1\sum_{j=1}^{l}j\int_0^{\alpha_{1,\xi}(j)}Q(Q\wedge a)^2\,\mathrm{d}\lambda+|\mathbb{E}(X_0X_0(a))|\sum_{j=1}^{l}j\int_0^{\alpha_{1,\xi}(j)}Q^2\,\mathrm{d}\lambda.$$

Finally, we have the simple bound $A_8\ll\int_0^1 Q^2\mathbb{1}_{Q>a}\,\mathrm{d}\lambda$.

### *4.5. End of the proof of Theorems 3.1 and 3.2*

We start with two preliminary results.

**Proposition 4.2.** *Let $(X_i)_{i\in\mathbb{Z}}$ be a stationary sequence of centered random variables in $\mathbb{L}^3$. Assume that the series $\sigma^2=\mathbb{E}(X_0^2)+2\sum_{k=1}^{\infty}\mathbb{E}(X_0X_k)$ converges absolutely and that $\sigma>0$. Consider the assumption* **H***: there exist positive constants $K$ and $M$ and a double array $(Y_{k,n})_{1\le k\le n}$ of independent and centered r.v.'s with common variance $\sigma^2$, such that, setting $T_n=Y_{1,n}+\cdots+Y_{n,n}$,*

$$d_1(S_n,T_n)\le M\quad\textit{and}\quad\max_{1\le k\le n}\mathbb{E}(|Y_{n,k}^3|)\le K^3.$$

*If* **H** *holds, then $d_1(S_n,\sqrt{n}\sigma Y)\le C$ for some constant $C$ depending only on $M$, $\sigma$ and $K$.*

**Proof.** Assume that **H** holds. Applying Theorem 5.17 in [23] we infer that there exists a constant $A$ such that, for any $x$,

$$|\mathbb{P}(T_n\le x\sqrt{n}\sigma)-\mathbb{P}(Y\le x)|\le A\left(\frac{K}{\sigma}\right)^3 n^{-1/2}(1+|x|)^{-3}.$$

Hence, integrating on the real line, $d_1(T_n,\sqrt{n}\sigma Y)\le AK^3\sigma^{-2}$. The result follows by taking $C=M+A\sigma^{-2}K^3$. □

We also need the following lemma, whose proof is elementary.

**Lemma 4.1.** *Let $\beta_2>0$ and $\beta_3$ be two fixed real numbers, and define*

$$m=\frac{\beta_3+\sqrt{\beta_3^2+\beta_2^3/2}}{\beta_2},\qquad m'=\frac{-\beta_2}{2m}\quad\textit{and}\quad t=\frac{\beta_2^3}{2\beta_2^3+4\beta_3(\beta_3+\sqrt{\beta_3^2+\beta_2^3/2})}.$$

*Let $Z_{\beta_2}$ and $B_{\beta_2,\beta_3}$ be two independent r.v.'s such that $Z_{\beta_2}$ has the distribution $\mathcal{N}(0,\beta_2/2)$ and $B_{\beta_2,\beta_3}$ is such that $\mathbb{P}(B_{\beta_2,\beta_3}=m)=t$ and $\mathbb{P}(B_{\beta_2,\beta_3}=m')=1-t$. Let $G_{\beta_2,\beta_3}=Z_{\beta_2}+B_{\beta_2,\beta_3}$. Then $\mathbb{E}(G_{\beta_2,\beta_3})=0$, $\mathbb{E}(G_{\beta_2,\beta_3}^2)=\beta_2$ and $\mathbb{E}(G_{\beta_2,\beta_3}^3)=\beta_3$.*

To prove Theorems 3.1 or 3.2, it is enough to see that under the assumptions of Theorems 3.1 or 3.2, the condition **H** of Proposition 4.2 holds. Without loss of generality, we assume that $\sigma^2=\mathbb{E}(X_0^2)+2\sum_{k=1}^{\infty}\mathbb{E}(X_0X_k)=1$ (the general case follows by dividing the random variables by $\sigma$). Denote by $b(l,a)$ the quantity $b(l)$ defined in (4.4) with $X_{k,1}=X_k(a)$ (see Section 4.3 for the definition of $X_k(a)$), and denote by $b(l,\infty)$ the quantity $b(l)$ with $X_{k,1}=X_k(\infty)=X_k$. Let $Y_{1,n},\dots,Y_{n,n}$ be $n$ independent random variables, independent of $(X_k)_{k\in\mathbb{Z}}$, such that $Y_{k,n}$ has the law of $G_{1,b(l(n,k),a(n,k))}$, where $G_{\beta_2,\beta_3}$ is defined in Lemma 4.1. Let $Y$ be a $\mathcal{N}(0,1)$-distributed random variable, independent of $(X_i, Y_{j,n})_{i\in\mathbb{Z},1\le j\le n}$, and let $T_n = Y_{1,n} + \cdots + Y_{n,n}$. Starting from (2.3), and keeping the same notations as in Notation 2.1, we have, as in Section 4.2,

$$\mathbb{E}(f(S_n) - f(T_n)) \le 2\mathbb{E}|Y| + \sum_{k=1}^{n} \mathbb{E}(\Delta_k). \tag{4.26}$$

By Lemma 2.1 applied with $B = B_{1,b(l(n,k+1),a(n,k+1))} + \cdots + B_{1,b(l(n,n),a(n,n))}$ we get that

$$\|f_k^{(i)}\|_\infty \le D_i(n-k+1)^{(1-i)/2}.$$

Define $\alpha^{-1}(u) = \sum_{i\ge 0} \mathbb{1}_{u<\alpha_{3,\boldsymbol{\xi}}(i)}$, and $R(u) = \alpha^{-1}(u)Q(u)$. Let $x_k = R^{-1}(\sqrt{k})$ and choose the truncation level $a(n,k) = \infty$ for Theorems 3.1(a), 3.2(a) and $a(n,k) = Q(x_{n-k+1})$ for Theorems 3.1(b), 3.2(b). Let $B_n$ be the set of positive integers $k$ such that $k-1 \le \sqrt{n-k+1}$ for Theorems 3.1(a) and 3.2(a), and $B_n$ be the set of positive integers $k$ such that $k-1 \le 4\alpha^{-1}(x_{n-k+1})$ for Theorems 3.1(b), 3.2(b). If $k$ belongs to $B_n$, take $l(n,k) = k-1$. If $k$ does not belong to $B_n$, take $l(n,k) = [\sqrt{n-k+1}]$ for Theorems 3.1(a), 3.2(a) and $l(n,k) = 4\alpha^{-1}(x_{n-k+1})$ for Theorems 3.1(b) and 3.2(b). Let $g(n) = \sup B_n$. Applying Proposition 4.1, with $Z = Y_{k,n}$, $\beta_2 = \sigma^2 = 1$, $\beta_3 = b(l(n,k), a(n,k))$ and $\beta_4 = \mathbb{E}(Y_{k,n}^4)$, we obtain that

$$\sum_{k=g(n)}^{n} \mathbb{E}(\Delta_k) \ll \sum_{k=g(n)}^{n} \sum_{i=1}^{9} \kappa_{i,n-k} A_i(k, l(n,k))$$

and

$$\sum_{k=1}^{g(n)-1} \mathbb{E}(\Delta_k) \ll \sum_{k=1}^{g(n)-1} \sum_{i=1}^{9} \kappa_{i,n-k} A_i(k, k-1),$$

where the numbers $\kappa_{i,m}$ are defined by $\kappa_{i,m} = m^{-i/2}$ for $i = 1,2,3,4$, $\kappa_{8,m} = \kappa_{2,m}$, $\kappa_{9,m} = \kappa_{3,m}$ and $\kappa_{i,m} = \kappa_{4,m}$ for $i = 5,6,7$. We only control the first term, the second one being easier to handle, since in that case $l(n,k) = k-1 \le \sqrt{n-k+1}$ for Theorems 3.1(a), 3.2(a) and $l(n,k) = k-1 \le 4\alpha^{-1}(x_{n-k+1})$ for Theorems 3.1(b), 3.2(b). To prove that condition **H** of Proposition 4.2 holds, it is enough to prove that for any $i$ in $[1,9]$,

$$\sup_{n>0} \sum_{k=g(n)}^{n} \kappa_{i,n-k} A_i(k, l(n,k)) < \infty. \tag{4.27}$$

The proof of (4.27) will be done using the upper bounds given in Section 4.3. We first prove Theorem 3.1 under condition (a) and next Theorem 3.1 under (b).

**Proof of Theorem 3.1 under condition (a).** In that case $a(n,k) = \infty$ and $l(n,k) = [\sqrt{n-k+1}]$. Consequently $A_8(k, l(n,k)) = A_9(k, l(n,k)) = 0$, so that we have to prove that (4.27) holds for $i$ in $[1,7]$.

By definition of $l(n,k)$, (4.27) holds for $i=1$ as soon as $\sum_{k=1}^{\infty} \theta_{0,1}([\sqrt{k}]) < \infty$, which is equivalent to condition (a) with $(p,q) = (0,1)$.

From (4.8), (4.27) holds for $i=2$ as soon as,

$$\text{for } (p,q) = (0,2) \text{ or } (1,2), \quad \sum_{k=1}^{\infty} \frac{1}{\sqrt{k}} \sum_{j=[\sqrt{k/2}]}^{\infty} \theta_{p,q}(j) < \infty. \tag{4.28}$$

Again these conditions are implied by condition (a) with $q=2$ and $p = 0,1$.

From (4.13), (4.27) holds for $i=3$ as soon as

$$\text{for } p=0,1,2 \text{ and } q=3, \quad \sum_{k=1}^{\infty} \frac{1}{\sqrt{k}} \sum_{j=[\sqrt{k/4}]}^{\infty} \theta_{p,q}(j) < \infty, \tag{4.29}$$

which holds true under condition (a) with $q=3$ and $p=0,1,2$.

From (4.7), (4.10), (4.19), (4.21) and (4.22) we infer that (4.27) holds for $i=5,6,7$ as soon as,

$$\text{for } q=4 \text{ and } 1 \le p \le 3, \quad \sum_{k=1}^{\infty} \frac{1}{k^{3/2}} \sum_{j=1}^{\infty} (2j \wedge \sqrt{k})^2 \theta_{p,q}(j) < \infty. \tag{4.30}$$

Clearly (4.30) holds as soon as

$$\sum_{k=1}^{\infty} \frac{1}{k^{3/2}} \sum_{j=1}^{[\sqrt{k}/2]} j^2 \theta_{p,q}(j) < \infty \quad \text{and} \quad \sum_{k=1}^{\infty} \frac{1}{\sqrt{k}} \sum_{j=[\sqrt{k}/2]}^{\infty} \theta_{p,q}(j) < \infty. \tag{4.31}$$

Interchanging the sums, we see that (4.31) holds under condition (a) with $q=4$ and $1 \le p \le 3$.

It remains to prove that (4.27) holds for $i=4$. Since $\sum_i \kappa_{i,4} \le \pi^2/6$, we infer that (4.27) holds for $i=4$ as soon as

$$\sup\{\mathbb{E}(Y_{k,n}^4)\colon\ 1 \le k \le n < \infty\} < \infty. \tag{4.32}$$

Now $Y_{k,n}$ and $Z_1 + B_{1,b(l(n,k),a(n,k))}$ have the same distribution, and consequently

$$\mathbb{E}(Y_{k,n}^4) \le 16(\mathbb{E}(Z_1^4) + \|B_{1,b(l(n,k),a(n,k))}\|_\infty^4). \tag{4.33}$$

Note that $\|B_{1,b(l(n,k),a(n,k))}\|_\infty = m \vee |m'|$, where $m$ and $m'$ have been defined in Lemma 4.1 with $\beta_2 = \sigma^2 = 1$ and $\beta_3 = b(l(n,k),a(n,k))$. Next $b(l(n,k),a(n,k)) \le b_3 < \infty$ with

$$b_3 = \theta_{0,3}(0) + 6\sum_{k=1}^{\infty} \theta_{1,3}(k) + 6\sum_{k=1}^{\infty} k\theta_{1,3}(k) + 6\sum_{k=1}^{\infty} k\theta_{2,3}(k). \tag{4.34}$$

Now $b_3$ is finite under condition (a). Hence, from Lemma 4.1

$$\|B_{1,b(l(n,k),a(n,k))}\|_\infty \le b_3 + \sqrt{b_3^2 + \frac{1}{2}}, \tag{4.35}$$

which completes the proof of (4.27) for $i=4$. □

**Proof of Theorem 3.1 under condition (b).** In that case we choose $a(n,k) = Q(x_{n-k+1})$ and $l(n,k) = 4\alpha^{-1}(x_{n-k+1})$. Note that, for any nonnegative measurable function $h$ and any positive $p$,

$$\sum_{i=0}^{\infty} i^{p-1} \int_0^{\alpha_{3,\boldsymbol{\xi}}(i)} h \, d\lambda < \infty \quad \text{if and only if} \quad \int_0^1 (\alpha^{-1})^p h \, d\lambda < \infty. \tag{4.36}$$

Hence (4.27) holds for $i=1$ as soon as,

$$\sum_{k=1}^{\infty} \int_0^1 \mathbb{1}_{k \le R^2} Q \, d\lambda < \infty,$$

which holds under (b).

Now recall that, for $i=2,3,6,7,9$, the terms $A_i(k,l(n,k))$ are decomposed into a sum of two terms: $A_i(k,l(n,k))=A_{i,1}(k,l(n,k))+A_{i,2}(k,l(n,k))$. Consequently, in order to prove that (4.27) holds for these values of $i$, we will prove that, for $j=1$ and $j=2$,

$$\sup_{n>0}\sum_{k=g(n)}^{n}\kappa_{i,n-k}A_{i,j}(k,l(n,k))<\infty. \tag{4.37}$$

For $i=2$, from (4.5), (4.37) holds for $i=2$ and $j=2$ as soon as

$$\int_0^1 \alpha^{-1}Q^2\sum_{k=1}^{\infty}\frac{1}{\sqrt{k}}\mathbb{1}_{k\le R^2}\,\mathrm{d}\lambda<\infty, \tag{4.38}$$

which follows from (b). From (4.8) and (4.9), (4.37) holds for $i=2$ and $j=1$ as soon as (4.38) holds. Hence (4.27) holds for $i=2$.

For $i=3$ and $i=9$, from (4.5), (4.37) holds for $j=2$ as soon as

$$\|X_0\|_1\sum_{k=1}^{\infty}\frac{1}{\sqrt{k}Q(x_k)}\sum_{j=1}^{\infty}\int_0^{\alpha_{3,\boldsymbol{\xi}}(j)}Q^2\mathbb{1}_{k\le R^2}\,\mathrm{d}\lambda<\infty,$$

which can be handled as (4.38) by noting that $\sqrt{k}Q(x_k)\ge\sqrt{k}Q(x_1)$. From (4.18) and (4.12), (4.37) holds for $i=3,9$ and $j=1$ as soon as

$$\sum_{k=1}^{\infty}\frac{R(x_k)}{k}\sum_{j=1}^{\infty}\int_0^{\alpha_{3,\boldsymbol{\xi}}(j)}Q^2\mathbb{1}_{k\le R^2}\,\mathrm{d}\lambda<\infty. \tag{4.39}$$

Since $R(x_k)\le\sqrt{k}$, (4.39) can be handled as (4.38). Hence (4.27) holds for $i=3$ and $i=9$.

For $i=6,7$ and from (4.5), (4.37) holds for $j=2$ as soon as

$$\|X_0\|_2^2\sum_{k=1}^{\infty}\frac{1}{\sqrt{k}Q^2(x_k)}\sum_{j=1}^{\infty}\int_0^{\alpha_{3,\boldsymbol{\xi}}(j)}Q^2\mathbb{1}_{k\le R^2}\,\mathrm{d}\lambda<\infty,$$

which can be handled as (4.38) by noting that $\sqrt{k}Q^2(x_k)\ge\sqrt{k}Q^2(x_1)$. From (4.10), (4.20), (4.22) and (4.24), (4.37) holds for $i=6,7$ and $j=1$ as soon as,

$$\sum_{k=1}^{\infty}\frac{1}{k^{3/2}}\sum_{j=1}^{\infty}j^2\mathbb{1}_{x_k\le\alpha_{3,\boldsymbol{\xi}}(j/4)}\int_0^{\alpha_{3,\boldsymbol{\xi}}(j)}Q^2(Q\wedge Q(x_k))^2\,\mathrm{d}\lambda<\infty. \tag{4.40}$$

Interchanging the sums and the integral, (4.40) holds as soon as

$$\int_0^1 Q^2\left(\sum_{k=1}^{\infty}\frac{(Q(x_k))^2}{k^{3/2}}\mathbb{1}_{k\le R^2}\sum_{j=1}^{\infty}j^2\mathbb{1}_{x_k\le\alpha_{3,\boldsymbol{\xi}}(j/4)}\right)\mathrm{d}\lambda<\infty, \tag{4.41}$$

$$\int_0^1(\alpha^{-1})^3Q^4\left(\sum_{k=1}^{\infty}\frac{1}{k^{3/2}}\mathbb{1}_{k>R^2}\right)\mathrm{d}\lambda<\infty. \tag{4.42}$$

Now (4.41) is equivalent to (4.38), and by definition of $R$, (4.42) follows from condition (b). Hence (4.27) holds for $i=6$ and $i=7$. In a similar way, for $i=5$, (4.27) can be derived from inequality (4.7).

We now prove (4.27) for $i=4$. First note that

$$\sum_{k=1}^{\infty}\frac{1}{k^{3/2}}|\mathbb{E}(X_0X_0^3(Q(x_k)))|\le\int_0^1 Q^2\sum_{k=1}^{\infty}\frac{Q^2(x_k)}{k^{3/2}}\mathbb{1}_{k\le R^2}\,\mathrm{d}\lambda+\int_0^1 Q^4\sum_{k>R^2}\frac{1}{k^{3/2}}\,\mathrm{d}\lambda$$

and these sums can be handled as in (4.41), (4.42). Since $\sum_i \kappa_{i,4} < \pi^2/6$, (4.27) holds for $i = 4$ as soon as (4.32) holds. Now as in the proof of (4.27) for $i = 4$ under condition (a), (4.33) holds and

$$b(l(n,k), a(n,k)) \le \sum_{k=1}^{\infty} 8k \int_0^{\alpha_{3,\xi}(k)} Q^3 \, d\lambda < \infty \tag{4.43}$$

under condition (b). Hence, from Lemma 4.1, (4.32) holds, which completes the proof of (4.27) holds for $i = 4$.

From (4.25), (4.27) holds for $A_8$ as soon as (4.38) holds. Hence Theorem 3.1 holds under (b). □

**Proof of Theorem 3.2.** Since the proof of Theorem 3.2(b) is similar to that of Theorem 3.1(b), we shall only give some hints at the end of this section.

To prove Theorem 3.2(a), we use the control of the $A_i$'s given in Section 4.4. Recall that, in that case, $A_1, A_5, A_8$ and $A_9$ are equal to zero. For $i = 2$, (4.27) holds as soon as

$$\sum_{k=1}^{\infty} \frac{1}{\sqrt{k}} \|E_0(X^2_{[\sqrt{k}]}) - 1\|_1 < \infty,$$

which follows from the first condition in (3.3). For $i = 6$, (4.27) holds as soon as

$$\sum_{k=1}^{\infty} \frac{1}{k\sqrt{k}} \sum_{j=1}^{\sqrt{k}} \|X_0^2(E_0(X_j^2) - 1)\|_1 < \infty,$$

which follows from the first condition in (3.3) by interchanging the sums. For $i = 3$, (4.27) holds as soon as

$$\sum_{k=1}^{\infty} \frac{1}{k} \sum_{j=[\sqrt{k}/2]}^{[\sqrt{k}]} \|X_0(E_0(X_j^2) - 1)\|_1 < \infty, \tag{4.44}$$

$$\sum_{k=1}^{\infty} \frac{1}{k} \sum_{j=[\sqrt{k}/2]}^{[\sqrt{k}]} \|E_0(X_j X^2_{[\sqrt{k}]}) - \mathbb{E}(X_j X^2_{[\sqrt{k}]})\|_1 < \infty. \tag{4.45}$$

Equation (4.44) follows from (3.3) by interchanging the sums. Equation (4.45) is equivalent to

$$\sum_{k=1}^{\infty} \frac{1}{k} \sum_{j=[k/2]}^{k} \|E_0(X_j X_k^2) - \mathbb{E}(X_j X_k^2)\|_1 < \infty,$$

which follows from the second condition in (3.3). For $i = 7$, (4.27) holds as soon as

$$\sum_{k=1}^{\infty} \frac{\|X_0\|_1}{k^{3/2}} \sum_{j=1}^{[\sqrt{k}]} j \|X_0(E_0(X_j^2) - 1)\|_1 < \infty, \tag{4.46}$$

$$\sum_{k=1}^{\infty} \frac{1}{k^{3/2}} \sum_{j=1}^{[\sqrt{k}]} \sum_{p=[j/2]}^{j} \|X_0(E_0(X_p X_j^2) - \mathbb{E}(X_p X_j^2))\|_1 < \infty, \tag{4.47}$$

$$\sum_{k=1}^{\infty} \frac{1}{k^{3/2}} \sum_{j=1}^{[\sqrt{k}]} \sum_{p=1}^{j} \|X_{-p}(E_0(X_j^2) - 1)\|_1 < \infty. \tag{4.48}$$

Interchanging the sums, we see that (4.46) and (4.48) follow from the first condition in (3.3), and (4.47) follows from the second condition in (3.3). For $i=4$, we proceed as in Theorem 3.1. We have the upper bound $b(l(n,k),\infty)\le d_3$ with,

$$d_3=\mathbb{E}(|X_0|^3)+3\sum_{k=1}^{\infty}\|X_0(E_0(X_k^2)-1)\|_1.$$

Hence (4.35) holds with $a(n,k)=\infty$, and the proof of $\mathbf{H}$ under (a) is complete.

The proof of Theorem 3.2(b) is similar to that of Theorem 3.1(b). For $1\le i\le 9$, (4.27) holds as soon as

$$\sum_{k=1}^{\infty}\frac{1}{\sqrt{k}}\int_0^1 Q^2\mathbb{1}_{k\le R^2}\,\mathrm{d}\lambda<\infty,\qquad \sum_{k=1}^{\infty}\frac{Q(x_k)}{k}\sum_{j=1}^{\infty}\mathbb{1}_{x_k\le\alpha(j/4)}\int_0^1 Q^2\mathbb{1}_{k\le R^2}\,\mathrm{d}\lambda<\infty$$

and

$$\sum_{k=1}^{\infty}\frac{1}{k^{3/2}}\sum_{j=1}^{\infty}j\mathbb{1}_{x_k\le\alpha(j/4)}\int_0^{\alpha(j)}Q^2(Q\wedge Q(x_k))^2\,\mathrm{d}\lambda<\infty.$$

Arguing as in the proof of Theorem 3.1(b), these inequalities follow from (3.4). □

## 5. Examples

### *5.1. Aperiodic Harris recurrent Markov chains*

Throughout this section, $K$ is a positive Harris recurrent Markov kernel on some separable state space $(E,\mathcal{E})$, i.e. there exists a unique probability measure $\pi$ with $\pi K=\pi$, and $K$ is $\pi$-recurrent. As in [5], $K$ is assumed to be aperiodic, which ensures that the stationary chain $(\xi_i)_{i\in\mathbb{Z}}$ with kernel $K$ is strongly mixing in the sense of Rosenblatt. Moreover, in the case of discrete Markov chains or chains with an atom, the rates of strong mixing and the integrability properties of the recurrence times are strongly linked, as proved by Theorem 2 in [3]: for any $r>-1$, $\sum_{k>0}k^r\alpha(k)<\infty$ if and only if $E(\tau^{r+2})<\infty$, where $\tau$ is the recurrence time (starting from the atom). From [26] the above series is convergent if and only if $\alpha^{-1}(u)=\sum_{i\ge 0}\mathbb{1}_{u<\alpha(i)}$ belongs to $L^{r+1}([0,1])$.

For any measurable function $f$, let $S_n(f)=f(\xi_1)+f(\xi_2)+\cdots+f(\xi_n)$. From Bolthausen's results ([3], Corollary 3 and [5], Theorem 1), the convergence rates in the Berry–Esseen theorem are $\mathrm{O}(n^{-1/2})$ as soon as

$$\pi(|f|^{3p})<\infty\quad\text{and}\quad\sum_{k>0}k^{(p+1)/(p-1)}\alpha(k)<\infty,\tag{5.1}$$

for any $p$ in $]1,\infty]$, provided that

$$\sigma^2=\pi(f^2)+2\sum_{n>0}\pi(fK^nf)>0.\tag{5.2}$$

From Theorem 3.1(b) above, we obtain the bound

$$d_1(n^{-1/2}S_n(f),\sigma Y)\le Cn^{-1/2}\tag{5.3}$$

as soon as $f$ satisfies (5.2) and (1.8), with $X_0=f(\xi_0)$ and $b=0$. From (4.36), the latter condition is equivalent to

$$\int_0^1[\alpha^{-1}(u)]^2Q^3_{|f(\xi_0)|}(u)\,\mathrm{d}u<\infty.\tag{5.4}$$

From the Hölder inequality applied with $s=p/(p-1)$ and $t=p$, we see that (5.4) holds as soon as (5.1) holds.

*Martingale difference sequences*
If $K(f)=0$ almost everywhere, the sequence $X_i = f(\xi_i)$ is a martingale difference sequence. Consequently Theorems 3.2 and 2.1 apply with $\sigma^2 = \pi(f^2)$. From Theorem 3.2(b), (5.3) holds as soon as the strong mixing coefficients satisfy (1.8) with $b=0$. Under the weaker condition

$$\int_0^{\alpha(k)} Q^3_{|f(\xi_0)|}(u)\,\mathrm{d}u = \mathrm{O}(k^{-\delta}), \tag{5.5}$$

Theorem 2.1(a) provides the rate

$$d_1(n^{-1/2}S_n(f), \sigma Y) = \mathrm{O}(n^{-\delta/2}). \tag{5.6}$$

When $f$ is a bounded function with $K(f)=0$ almost everywhere, (5.3) holds under the summability condition $\sum_k \alpha(k) < \infty$, which is related to the ergodicity of degree 2 (cf. [21], Section 6.4). From (5.5), the rate (5.6) holds under the weaker condition $\alpha(k) = \mathrm{O}(k^{-\delta})$.

### 5.2. *The transformation $\Theta(x) = 2x - [2x]$*

Let $\lambda$ be the Lebesgue measure on $[0,1]$ and consider the map $\Theta$ from $[0,1]$ to $[0,1]$: $\Theta(x) = 2x - [2x]$. On the probability space $([0,1], \lambda)$, the sequence $(\Theta^i)_{i>0}$ is strictly stationary. Note also that $(\Theta, \Theta^2, \dots, \Theta^n)$ is distributed as $(\xi_n, \dots, \xi_1)$, where $(\xi_i)_{i\in\mathbb{Z}}$ is a Markov chain with invariant distribution $\lambda$ and transition kernel

$$Kf(x) = \frac{1}{2}\left(f\left(\frac{x}{2}\right) + f\left(\frac{x+1}{2}\right)\right).$$

Hence, we can obtain information on the distribution of $S_n(f) = f\circ\Theta + \cdots + f\circ\Theta^n$ by studying that of $f(\xi_1) + \cdots + f(\xi_n)$. For instance, we can apply the criterion of Dedecker and Rio [8] for the central limit theorem: if $\lambda(f) = 0$,

$$\lambda(f^2) < \infty \quad \text{and} \quad \sum_{k>0} \lambda(|fK^k(f)|) < \infty, \tag{5.7}$$

then $\sigma^2 = \lambda(f^2) + 2\sum_{k=1}^{\infty} \lambda(f\cdot f\circ\Theta^k)$ converges absolutely, and $n^{-1/2}S_n(f)$ converges in distribution to a Gaussian random variable with mean 0 and variance $\sigma^2$. Now it is easy to see that (5.7) holds as soon as, for some $p \in [2,\infty]$,

$$f \in \mathbb{L}^p(\lambda) \quad \text{and} \quad \int_0^1 \frac{1}{t} w_{p/(p-1)}(t)\,\mathrm{d}t < \infty, \tag{5.8}$$

where $w_q(f,t)$ is the $\mathbb{L}^q([0,1],\lambda)$-modulus of continuity of $f$ in $\mathbb{L}^q([0,1],\lambda)$. For $p=2$, the criterion (5.8) has been obtained by Ibragimov [16]. For $p=\infty$, the criterion (5.8) follows from the $\mathbb{L}^1$-criterion of Gordin [14] applied to sequences of bounded variables.

In the same way, applying Theorems 2.2, 3.2 and 3.1, we obtain the following result:

**Theorem 5.1.** *Let $f$ be a measurable function from $[0,1]$ to $\mathbb{R}$ such that $\lambda(f) = 0$.*

(a) *Assume that, for some $p \in [3,\infty]$,*

$$f \in \mathbb{L}^p(\lambda) \quad \textit{and} \quad \int_0^1 \frac{|\log t|}{t} w_{p/(p-2)}(f,t)\,\mathrm{d}t < \infty. \tag{5.9}$$

*If $\sigma > 0$, there exists a constant $C$ such that, $d_1(S_n(f), \sqrt{n}\sigma Y) \le C\log n$.*

(b) *Assume that, for some* $p \in [4, \infty]$,

$$f \in \mathbb{L}^p(\lambda) \quad \text{and} \quad \int_0^1 \frac{|\log t|}{t} w_{p/(p-3)}(f,t)\,\mathrm{d}t < \infty. \tag{5.10}$$

*If* $\sigma > 0$, *then there exists a constant* $C$ *such that,*

$$d_1(S_n(f), \sqrt{n}\sigma Y) \le C. \tag{5.11}$$

(c) *Assume that* $f(x + (1/2)) = -f(x)$ *for almost every* $x \in [0, 1/2]$. *Then the sequence* $(f(\xi_n))_{n\in\mathbb{Z}}$ *is a stationary martingale difference sequence, so that* $\sigma^2(f) = \lambda(f^2)$. *If moreover, for some* $p \in [4, \infty]$,

$$\lambda(|f|^p) < \infty \quad \text{and} \quad \int_0^1 \frac{1}{t} w_{p/(p-3)}(f,t)\,\mathrm{d}t < \infty,$$

*then (5.11) holds.*

(d) *Assume that* $f = f_1 - f_2$, *where* $f_1$ *and* $f_2$ *are nondecreasing functions. Assume moreover that*

$$\int_0^1 (\log(t-t^2))^2 |f_1(t)|^3\,\mathrm{d}t < \infty \quad \text{and} \quad \int_0^1 (\log(t-t^2))^2 |f_2(t)|^3\,\mathrm{d}t < \infty.$$

*If* $\sigma > 0$, *then (5.11) holds.*

***Remark 5.1.*** *If* $f$ *belongs to* $\mathbb{L}^3(\lambda)$, *Ibragimov [17] obtained the Berry–Esseen type estimate*

$$\sup_{x\in\mathbb{R}} |\mathbb{P}(S_n(f) \le x\sqrt{n}\sigma) - \mathbb{P}(Y \le x)| \le C\left(\frac{\log n}{n}\right)^{1/2}, \tag{5.12}$$

*under the condition* $w_3(f,t) \le Ct^\alpha$ *for some* $\alpha > 0$. *This condition is slightly stronger than our condition (5.9) with* $p = 3$. *Applying Theorem 9 in [19], one can obtain the bound* $Cn^{-1/2}$ *in (5.12) as soon as*

$$f \in \mathbb{L}^\infty(\lambda) \quad \text{and} \quad \int_0^1 \frac{|\log t|}{t} w_\infty(f,t)\,\mathrm{d}t < \infty,$$

*where* $w_\infty(f,t)$ *is the modulus of continuity of* $f$.

***Remark 5.2.*** *If* $f \in \mathbb{L}^\infty(\lambda)$ *and if* $K(f) = 0$ *almost everywhere, then (5.11) holds under the criterion (5.8) applied to* $p = \infty$.

***Remark 5.3.*** *Applying the Hausdorff–Young inequality (cf. [15], p. 202), Ibragimov [16, 17] proved that (5.8) holds for* $p = 2$ *as soon as the Fourier coefficients of* $f$ *satisfy* $|\hat{f}(n)| \le Mn^{-1/2}(\log(n))^{-3/2-\epsilon}$ *for some positive* $M$ *and* $\epsilon$, *and that (5.12) holds as soon as* $|\hat{f}(n)| \le Mn^{-2/3-\epsilon}$. *Using the same arguments, one can prove that (5.9) holds for* $p = 3$ *as soon as* $|\hat{f}(n)| \le Mn^{-2/3}(\log(n))^{-8/3-\epsilon}$, *and that (5.10) holds for* $p = 4$ *as soon as* $|\hat{f}(n)| \le Mn^{-3/4}(\log(n))^{-11/4-\epsilon}$.

**Proof of Theorem 5.1.** Point (a) follows from Theorem 2.2 and point (b) follows from Theorem 3.1(a). The proofs being similar, we shall only prove point (b). Let us just see how to control the coefficient $\theta_{1,4}(l)$, the other one being easier to handle. The sequence $(\xi_i)_{i\in\mathbb{Z}}$ being a stationary Markov chain with invariant distribution $\lambda$ and transition kernel $K$, the coefficient $\theta_{1,4}(l)$ is equal to

$$\sup_{k\ge l, i\ge 0, j\ge 0} \int_0^1 \left| f(x)\left(K^k(fK^i(fK^j(f)))(x) - \int_0^1 (fK^i(fK^j(f)))(x)\,\mathrm{d}x\right)\right| \mathrm{d}x.$$

From Theorem 1 in [12], we infer that, for any $h$ in $\mathbb{L}^q([0,1],\lambda)$,

$$\left\|\int |K^k(h)(x)-\lambda(h)|^q\,\mathrm{d}x\right\|_{q,\lambda} \le 2w_q(h,2^{-k}).$$

Hence, applying Hölder's inequality, we obtain that

$$\theta_{1,4}(l) \le \sup_{i\ge 0, j\ge 0} 2\|f\|_{p,\lambda} w_{p/(p-1)}(fK^i(fK^j(f)),2^{-l}).$$

We now use the elementary facts that, for $p\ge q$ and $r\ge q$,

$$w_q(fg,t) \le \|f\|_{p,\lambda} w_{pq/(p-q)}(g,t) + \|g\|_{r,\lambda} w_{rq/(r-q)}(f,t),$$

and that $w_q(K(f),t)\le w_q(f,t)$. It follows that

$$\begin{aligned}\theta_{1,4}(l) &\le \sup_{j\ge 0} 2\|f\|_{p,\lambda}^2 (w_{p/(p-2)}(fK^j(f),2^{-l}) + \|f\|_{p,\lambda} w_{p/(p-3)}(f,2^{-l}))\\ &\le 6\|f\|_{p,\lambda}^3 w_{p/(p-3)}(f,2^{-l}).\end{aligned}$$

Hence, if $f$ belongs to $\mathbb{L}^p([0,1],\lambda)$ for some $p\ge 4$, $\sum_{l>0} l\theta_{1,4}(l)$ is finite as soon as $\sum_{l>0} l\, w_{p/(p-3)}(f,2^{-l})$ is finite, which is equivalent to the condition of (b).

To prove (c), note that, if $f(x+(1/2)) = -f(x)$ for almost every $x\in[0,1/2]$, then $(f(\xi_i))_{i\in\mathbb{Z}}$ is a sequence of martingale differences, so that Theorem 3.2(a) applies. To conclude, use the control of $\theta_{i,j}(l)$ given above.

It remains to prove (d). Let $BV_a$ be the space of left continuous bounded variation functions $f$ on $[0,1]$ such that $\|\mathrm{d}f\|_v\le a$ (here $\|\cdot\|_v$ is the variation norm). Let $f^{(0)} = f-\lambda(f)$ and $\mathcal{M}_0 = \sigma(Y_i, i\le 0)$. Arguing as in Lemma 1 of [7], one can see that, for any $i_l > \cdots > i_1 > n$,

$$\alpha(\mathcal{M}_0,(\xi_{i_1},\dots,\xi_{i_l})) = \sup_{f_1,\dots,f_l\in BV_1}\left\|\mathbb{E}\left(\prod_{j=1}^l f_j^{(0)}(\xi_{i_j})\Big|\mathcal{M}_0\right) - \mathbb{E}\left(\prod_{j=1}^l f_j^{(0)}(\xi_{i_j})\right)\right\|_1.$$

Since $K$ maps $BV_1$ to $BV_{1/2}$, we infer that $f^{(0)}\cdot(K^i(g))^{(0)}$ belongs to $BV_1$ for any $i>0$ and any $f,g$ in $BV_1$. It follows that, for any $i_l>\cdots>i_1\ge n$,

$$\alpha(\mathcal{M}_0,(\xi_{i_1},\dots,\xi_{i_l})) \le \alpha(\mathcal{M}_0,\xi_{i_1}) \le 2^{-n},$$

so that $\alpha_{3,\boldsymbol{\xi}}(n)\le 2^{-n}$. Applying Theorem 3.1(a), $d_1(S_n(f),\sqrt{n}\sigma Y)\le C$ as soon as

$$\int_0^1 (\log t)^2 Q_{|f_1|}^3(t)\,\mathrm{d}t < \infty \quad \text{and} \quad \int_0^1 (\log t)^2 Q_{|f_2|}^3(t)\,\mathrm{d}t<\infty,$$

where $Q_f$ is the generalized inverse of $t\to\lambda(f>t)$. Let $f^+ = f\vee 0$ and $f^- = -(f\wedge 0)$. By Lemma 2.1 in [26],

$$\int_0^1 (\log t)^2 Q_{|f_1|}^3(t)\,\mathrm{d}t \le \int_0^1 (\log t)^2 Q_{f_1^+}^3(t)\,\mathrm{d}t + \int_0^1 (\log t)^2 Q_{f_1^-}^3(t)\,\mathrm{d}t.$$

Clearly $Q_{f_1^+}(t) = f_1^+(1-t)$ almost everywhere and $Q_{f_1^-(t)} = f_1^-(t)$ almost everywhere. Of course the same is true with $f_2$ and the result follows. □

*5.3. Symmetric random walk on the circle*

Let $K$ be the Markov kernel defined by $Kf(x) = (f(x+a) + f(x-a))/2$ on $T = \mathbb{R}/\mathbb{Z}$, with $a$ irrational in $[0,1]$. The Lebesgue–Haar measure $m$ is invariant under $K$. Furthermore $K$ is a symmetric operator on $\mathbb{L}^2(m)$, and consequently the Kipnis–Varadhan or the Gordin–Lifshitz central limit theorems apply. Let $(\xi_i)_{i\in\mathbb{Z}}$ be the stationary Markov chain with transition kernel $K$. For $f$ in $\mathbb{L}^2(m)$ with $m(f) = 0$, set

$$S_n(f) = f(\xi_1) + f(\xi_2) + \cdots + f(\xi_n). \tag{5.13}$$

Then the central limit theorem holds for $n^{-1/2}S_n(f)$ as soon as the series of covariances

$$\sigma^2 = \int_T f^2 \,\mathrm{d}m + 2\sum_{n>0}\int_T fK^n f \,\mathrm{d}m \tag{5.14}$$

is convergent and the limiting distribution is $\mathcal{N}(0,\sigma^2)$ (cf. [9], Section 2). Our aim in this section is to give conditions on $f$ and on the properties of the irrational number $a$ ensuring optimal rates of convergence in the central limit theorem.

***Definition 5.1.*** *$a$ is said to be badly approximable by rationals if for any positive $\varepsilon$, the inequality $d(ka,\mathbb{Z}) < |k|^{-1-\varepsilon}$ has only finitely many solutions for $k \in \mathbb{Z}$.*

From Roth's theorem the algebraic numbers are badly approximable (cf. [27]). Note also that the set of badly approximable numbers in $[0,1]$ has Lebesgue measure 1. We will now give results for the symmetric random walk on the circle in the case of badly approximable numbers $a$.

**Theorem 5.2.** *Suppose that $a$ is badly approximable by rationals. Let $f$ be a function in $\mathbb{L}^2(m)$ with $m(f) = 0$ and $m(f^2) > 0$.*

(a) *If the Fourier coefficients $\hat f(k)$ of $f$ satisfy $\sup_{k\neq 0}|k|^{1+\varepsilon}|\hat f(k)| < \infty$ for some positive $\varepsilon$, then $n^{-1/2}S_n(f)$ converges in distribution to a nondegenerate Gaussian distribution $\mathcal{N}(0,\sigma^2)$.*

(b) *If the Fourier coefficients $\hat f(k)$ of $f$ satisfy $\sup_{k\neq 0}|k|^{4+\varepsilon}|\hat f(k)| < \infty$ for some positive $\varepsilon$, then*

$$\sup_{x\in\mathbb{R}}|\mathbb{P}(S_n \le x\sigma\sqrt{n}) - \mathbb{P}(Y \le x)| = \mathrm{O}(n^{-1/2}), \tag{5.15}$$

$$d_1(n^{-1/2}S_n, \sigma Y) = \mathrm{O}(n^{-1/2}). \tag{5.16}$$

***Remark 5.4.*** *The assumption $\hat f(k) = \mathrm{O}(|k|^{-1-\varepsilon})$ in Theorem 5.2*(a) *implies that $f$ is $\varepsilon$-Hölderian, and therefore uniformly continuous. Conversely, if $f$ is $C^{1+\varepsilon}$ then $f$ satisfies* (a). *In the same way the condition $\hat f(k) = \mathrm{O}(|k|^{-4-\varepsilon})$ in* (b) *implies that $f$ is $C^{3+\varepsilon}$ and conversely any $C^{4+\varepsilon}$ function $f$ satisfies* (b).

**Proof of Theorem 5.2.** Since

$$\int_T fK^n f \,\mathrm{d}m = \sum_{k\in\mathbb{Z}^*}\cos^n(2\pi ka)|\hat f(k)|^2,$$

the series in (5.14) is convergent if $\sum_{k\in\mathbb{Z}^*}\cot^2(\pi ka)|\hat f(k)|^2 < \infty$. Moreover, interverting the sums, we get that $\sigma^2 = \sum_{k\in\mathbb{Z}^*}\cot^2(\pi ka)|\hat f(k)|^2$. Since $\cot^2(\pi ka) > 0$ for any $k$ in $\mathbb{Z}^*$, it ensures that $\sigma^2 > 0$.

When $\{ka\} = d(ka,\mathbb{Z})$ tends to 0, $\cot^2(\pi ka) \sim \pi^{-2}\{ka\}^{-2}$, so that the convergence of the series in (5.14) is equivalent to

$$\sum_{k\in\mathbb{Z}^*}\{ka\}^{-2}|\hat f(k)|^2 < \infty, \tag{5.17}$$

as shown in [9].

In order to complete the proof of Theorem 5.2(a), we will need the elementary fact below.

**Lemma 5.1.** *Let $a$ be a badly approximable irrational number. Then, for any positive $\eta$, there exists some positive constant $C$ such that, for any nonnegative integer $N$ and any $p \geq 2$, $\sum_{k\in[2^N,2^{N+1}[}\{ka\}^{-p} \leq 2C^p 2^{p(N+2)(1+\eta)}$.*

**Proof.** Let $k$ and $l$ be integers in $I_N = [2^N, 2^{N+1}[$ with $k \neq l$. From the equality $|\{ka\} - \{la\}| = \min(\{(l-k)a\}, \{(l+k)a\})$ and Definition 5.1, we get that $|\{ka\} - \{la\}| \geq C^{-1}|k-l|^{-1-\eta} \geq C^{-1}2^{-(N+2)(1+\eta)}$ for some positive constant $C$. Now, denoting by $x_1^N, \dots, x_{2^N}^N$ the order statistic of $(\{ka\})_{k\in I_N}$,

$$x_m^N \geq x_1^N + (m-1)C^{-1}2^{-N(2+\eta)} \geq mC^{-1}2^{-(2+N)(1+\eta)}.$$

Hence

$$\sum_{k=2^N}^{2^{N+1}-1} \{ka\}^{-p} = \sum_{m=1}^{2^N} (x_m^N)^{-p} \leq C^p 2^{p(N+2)(1+\eta)} \sum_{m=1}^{2^N} m^{-p},$$

which implies Lemma 5.1. □

Now, applying Lemma 5.1 with $p = 2$ and $\eta = \varepsilon/2$, we get that

$$\sum_{k\in I_N} \{ka\}^{-2}(|\hat{f}(k)|^2 + |\hat{f}(-k)|^2) \leq 4C^2 2^{N(2+\varepsilon)} \max_{k\in I_N} |\hat{f}(k)|^2 \leq C' 2^{-N\varepsilon}$$

under the assumptions of Theorem 5.2(a), which implies the convergence of the series in (5.17). Therefore Theorem 5.2(a) holds.

We now prove Theorem 5.2(b). Equation (5.15) is a byproduct of Jan's theorem ([19], Theorem 9, page 61 or [20], Theorem 1) and (5.16) is a corollary of our estimates of the minimal $\mathbb{L}^1$-distance. The main tool is Lemma 5.2.

***Notation 5.1.*** *For $s > 0$, let $\mathcal{F}_s$ be the class of 1-periodic functions $g$ such that $\hat{g}(0) = 0$ and $|\hat{g}(k)| \leq |k|^{-s}$ for any $k$ in $\mathbb{Z}^*$.*

**Lemma 5.2.** *Let $a$ be a badly approximable irrational number. Then, for any $\varepsilon$ in $]0,1]$*

$$\sum_{n>0} n \sup_{g\in\mathcal{F}_{4+4\varepsilon}} \|K^n g\|_\infty < \infty.$$

**Proof.** For $g$ in $\mathbb{L}^2(m)$ with $m(g) = 0$,

$$K^n g(x) = \sum_{k\in\mathbb{Z}^*} \cos^n(2\pi ka)\hat{g}(k)\exp(2i\pi kx).$$

Therefore

$$\sup_{g\in\mathcal{F}_{4+4\varepsilon}} \|K^n g\|_\infty \leq \sum_{k\in\mathbb{Z}^*} |\cos^n(2\pi ka)||k|^{-4(1+\varepsilon)},$$

which ensures that

$$\sum_{n>0} n \sup_{g\in\mathcal{F}_{4+4\varepsilon}} \|K^n g\|_\infty \leq \sum_{k\in\mathbb{Z}^*} (1 - |\cos(2\pi ka)|)^{-2}|k|^{-4(1+\varepsilon)} \leq \sum_{k\in\mathbb{Z}^*} (|k|^{1+\varepsilon}\{2ka\})^{-4}.$$

Next, applying Lemma 5.1 with $\eta = \varepsilon/2$ and $p = 4$, we get that

$$\sum_{k\in\mathbb{Z}^*} (|k|^{1+\varepsilon}\{2ka\})^{-4} \le 4C^4 \sum_{N\ge 0} 2^{(4+2\varepsilon)(N+2)} \max_{k\in I_N} k^{-4(1+\varepsilon)} < \infty,$$

which implies Lemma 5.2. □

We now complete the proof of Theorem 5.2(b). Set $X_p = f(\xi_p)$. In view of the Berry–Esseen type Theorem 9 in [19] and Theorem 3.1 we have to bound up the coefficients

$$\psi_n = \sup\{\|E_0(X_{p_1}\cdots X_{p_j}) - \mathbb{E}(X_{p_1}\cdots X_{p_j})\|_\infty\colon j\le 3, n\le p_1\le\cdots\le p_j\}$$

in such a way that $\sum_n n\psi_n < \infty$.

We proceed as in [19]. Set $p_0 = n$. Then (in the case $j = 3$)

$$E_0(X_{p_1}\cdots X_{p_j}) = E_0(E_{p_0}(X_{p_1}(E_{p_1}(X_{p_2}E_{p_2}(X_{p_3}))))).$$

Hence, setting $q_i = p_i - p_{i-1}$, we get $E_0(X_{p_1}\cdots X_{p_j}) = K^n(K^{q_1}(fK^{q_2}(fK^{q_3}f)))$. Starting from this equality, we now prove that, for $s > 1$ there exists some constant $C_s$ (depending only on $s$) such that, for any $f\in\mathcal{F}_s$,

$$\psi_n \le C_s \sup_{g\in\mathcal{F}_s} \|K^n g\|_\infty. \tag{5.18}$$

To prove (5.18) one can prove that, for $f$ in $\mathcal{F}_s$ and $g = K^{q_1}(fK^{q_2}(f\cdots K^{q_j}f)\cdots)$,

$$|\hat{g}(k)| \le C_s|k|^{-s} \quad \text{for any } k\in\mathbb{Z}^*, \tag{5.19}$$

any $j\le 3$ and all natural integers $q_1,\ldots,q_j$. This is derived from Lemma 5.3.

**Lemma 5.3.** *Let $s > 1$. For any $g$ in $\mathcal{F}_s$ and any natural $p$, $K^p g$ lies in $\mathcal{F}_s$. For any $g$ and $h$ in $\mathcal{F}_s$ and any $k\ne 0$, $|\widehat{gh}(k)| \le (s-1)^{-1}2^{2s+1}|k|^{-s}$.*

The proof of Lemma 5.3, being elementary, is omitted. Now, from (5.18) and Lemma 5.2, $\sum_n n\psi_n < \infty$ under the assumptions of Theorem 5.2(b). Since the function $f$ is uniformly bounded, it implies (5.15) via Theorem 9 in [19] and (5.16) via Theorem 3.1. □

## Appendix

In this section, we give an upper bound for the expectation of the product of $k$ centered random variables $\prod_{i=1}^k (X_i - \mathbb{E}(X_i))$. This upper bound is given in term of a dependence coefficients $\alpha(X_1,\ldots,X_k)$, which is a generalization of the coefficient introduced in [26], Eq. (1.8a), for $k = 2$ (note that, for $k = 2$, our definition differs from that of Rio by a factor 2).

**Proposition A.1.** *Let $X = (X_1,\ldots,X_k)$ be a random variable with values in $\mathbb{R}^k$ and define the number*

$$\alpha = \alpha(X_1,\ldots,X_k) = \sup_{(x_1,\ldots,x_k)\in\mathbb{R}^k} \left|\mathbb{E}\left(\prod_{i=1}^k \mathbb{1}_{X_i > x_i} - \mathbb{P}(X_i > x_i)\right)\right|. \tag{A.1}$$

*Let $F_i$ be the distribution function of $X_i$, let $F_i^{-1}$ be the generalized inverse of $F_i$ and let $D_i(u) = (F_i^{-1}(1-u) - F_i^{-1}(u))_+$. We have the inequality*

$$\left|\mathbb{E}\left(\prod_{i=1}^k X_i - \mathbb{E}(X_i)\right)\right| \le 2\int_0^{\alpha/2} \left(\prod_{i=1}^k D_i(u)\right) du. \tag{A.2}$$

*In particular, if $X_1$ is $\mathcal{M}$-measurable, we have $\alpha \le \alpha(\mathcal{M},(X_2,\dots,X_k))$. Hence*

$$\left|\mathbb{E}\left(\prod_{i=1}^k X_i - \mathbb{E}(X_i)\right)\right| \le 2\int_0^{\alpha(\mathcal{M},(X_2,\dots,X_k))/2}\left(\prod_{i=1}^k D_i(u)\right)\mathrm{d}u. \tag{A.3}$$

**Proof.** We have that

$$\mathbb{E}\left(\prod_{i=1}^k X_i - \mathbb{E}(X_i)\right) = \int \mathbb{E}\left(\prod_{i=1}^k \mathbb{1}_{X_i>x_i} - \mathbb{P}(X_i>x_i)\right)\mathrm{d}x_1\cdots\mathrm{d}x_k. \tag{A.4}$$

Now $A = |\mathbb{E}(\prod_{i=1}^k(\mathbb{1}_{X_i>x_i} - \mathbb{P}(X_i>x_i)))|$ is such that $A \le \alpha$, and for any $1 \le i \le k$,

$$A \le 2\mathbb{P}(X_i>x_i)\mathbb{P}(X_i\le x_i)\wedge\alpha \le 2\left\{\mathbb{P}(X_i>x_i)\wedge\mathbb{P}(X_i\le x_i)\wedge\frac{\alpha}{2}\right\}. \tag{A.5}$$

Consequently, we obtain from (A.4) and (A.5) that

$$\begin{aligned}\left|\mathbb{E}\left(\prod_{i=1}^k X_i - \mathbb{E}(X_i)\right)\right| &\le 2\int_0^{\alpha/2}\left(\prod_{i=1}^k\int \mathbb{1}_{u<\mathbb{P}(X_i>x_i)}\mathbb{1}_{u\le\mathbb{P}(X_i\le x_i)}\,\mathrm{d}x_i\right)\mathrm{d}u\\ &\le 2\int_0^{\alpha/2}\left(\prod_{i=1}^k\int \mathbb{1}_{F_i^{-1}(u)\le x_i<F_i^{-1}(1-u)}\,\mathrm{d}x_i\right)\mathrm{d}u\end{aligned}$$

and (A.2) follows. □

**Lemma A.1.** *Let $X_+ = \max(0,X)$ and $X_- = -\min(0,X)$. For almost every $u<1/2$, we have the inequalities $0 \le D_X(u) \le Q_{X_+}(u) + Q_{X_-}(u) \le 2Q_{|X|}(u)$. Furthermore the second inequality is an equality if 0 is a median for $X$.*

**Proof.** First, we have $F_X^{-1}(1-u) = Q_X(u) \le Q_{X_+}(u)$. Next, by definition of $F_X^{-1}$, we have $-F_X^{-1}(u) = \sup\{x\colon\ \mathbb{P}(-X\ge x)\ge u\}$. By definition $Q_{-X}(u) = \inf\{x\colon\ \mathbb{P}(-X>x)\ge u\}$, so that $-F_X^{-1}(u) = Q_{-X}(u)$ for every continuity point $u$ of $Q_{-X}$ and hence almost everywhere. To obtain the desired inequality, note that $Q_{-X}(u) \le Q_{X_-}(u)$. □

**Corollary A.1.** *Let $X=(X_1,\dots,X_k)$ be a random variable with values in $\mathbb{R}^k$ and let $\alpha$ be defined by (A.1). Let $(f_i)_{1\le i\le k}$ be $k$ functions from $\mathbb{R}$ to $\mathbb{R}$, such that $f_i = f_i^{(1)} - f_i^{(2)}$ where $f_i^{(1)}$ and $f_i^{(2)}$ are nondecreasing. For $1\le i\le k$ and $j\in\{1,2\}$, let $Q_i^{(j)} = Q_{|f_i^{(j)}(X_i)|}$. We have the inequality*

$$\left|\mathbb{E}\left(\prod_{i=1}^k f_i(X_i) - \mathbb{E}(f_i(X_i))\right)\right| \le 2^{k+1}\sum_{j_1=1}^2\cdots\sum_{j_k=1}^2\int_0^{\alpha/2} Q_1^{(j_1)}(u)\cdots Q_k^{(j_k)}(u)\,\mathrm{d}u.$$

*In particular, if $X_1$ is $\mathcal{M}$-measurable,*

$$\left|\mathbb{E}\left(\prod_{i=1}^k f_i(X_i) - \mathbb{E}(f_i(X_i))\right)\right| \le 2^{k+1}\sum_{j_2=1}^2\cdots\sum_{j_k=1}^2\int_0^{\alpha(\mathcal{M},(X_2,\dots,X_k))/2} Q_{|f_1(X_1)|}(u)\left(\prod_{i=2}^k Q_i^{(j_i)}(u)\right)\mathrm{d}u.$$

**Proof.** Clearly

$$\left|\mathbb{E}\left(\prod_{i=1}^k f_i(X_i) - \mathbb{E}(f_i(X_i))\right)\right| \le \sum_{j_1=1}^2\cdots\sum_{j_k=1}^2\left|\mathbb{E}\left(\prod_{i=1}^k f_i^{(j_i)}(X_i) - \mathbb{E}(f_i^{(j_i)}(X_i))\right)\right|. \tag{A.6}$$

Since $f_i^{(j_i)}$ is nondecreasing, $\alpha(f_1^{(j_1)}(X_1),\ldots,f_k^{(j_k)}(X_k)) \leq \alpha(X_1,\ldots,X_k)$. To obtain the result, apply (A.2) and Lemma A.1 to each term of the sum in (A.6). □